\documentclass[12pt]{amsart}
\usepackage{amssymb,amscd,amsthm,verbatim}

\newtheorem{prop}{Proposition}
\newtheorem{thm}[prop]{Theorem}

\newtheorem{lem}[prop]{Lemma}
\newtheorem{ex}{Example}

\newtheorem{conj}{Conjecture}

\theoremstyle{definition}

\numberwithin{equation}{section}

\begin{document}

\title{A generalization of the Kostka-Foulkes polynomials}

\author{Anatol N. Kirillov}
\address{Steklov Mathematical Institute\\St.Petersburg\\Russia}
\email{kirillov@pdmi.ras.ru}

\author{Mark Shimozono}
\address{Department of Mathematics\\Virginia Tech\\Blacksburg, VA}
\email{mshimo@math.vt.edu}

\begin{abstract} Combinatorial objects called rigged
configurations give rise to $q$-analogues of certain
Littlewood-Richardson coefficients.
The Kostka-Foulkes polynomials and two-column
Macdonald-Kostka polynomials occur as special cases.
Conjecturally these polynomials coincide with the Poincar\'e polynomials
of isotypic components of certain graded $GL(n)$-modules supported
in a nilpotent conjugacy class closure in $gl(n)$.
\end{abstract}

\maketitle

\newcommand{\gl}{\mathrm{gl}}
\newcommand{\LRC}{LR}
\newcommand{\Rh}{\widehat R}
\newcommand{\cat}{\mathrm{Cat}}
\newcommand{\ccat}{\mathrm{CCat}}
\newcommand{\key}{\text{key}}
\newcommand{\charge}{\text{charge}}
\newcommand{\cocharge}{\text{cocharge}}
\newcommand{\emptytab}{\emptyset}
\newcommand{\la}{\lambda}
\newcommand{\Hom}{\mathrm{Hom}}
\newcommand{\Ind}{\mathrm{Ind}}
\newcommand{\Res}{\mathrm{Res}}
\newcommand{\Roots}{\mathrm{Roots}}
\newcommand{\word}{\mathrm{word}}
\newcommand{\rows}{\mathrm{rows}}
\newcommand{\cols}{\mathrm{columns}}
\newcommand{\shape}{\mathrm{shape}}
\newcommand{\ev}{\mathrm{ev}}
\newcommand{\rev}{\mathrm{rev}}
\newcommand{\Lmax}{L_{\mathrm{max}}}
\newcommand{\C}{\mathbb{C}}
\newcommand{\NN}{\mathbb{N}}
\newcommand{\Z}{\mathbb{Z}}
\newcommand{\dom}{\trianglerighteq}
\newcommand{\SR}{\mathcal{R}}
\newcommand{\qbinom}[2]{\binom{#1}{#2}_q}
\newcommand{\mh}{\widehat{m}}
\newcommand{\mt}{\widehat{t}}
\newcommand{\alh}{\widehat{\alpha}}
\newcommand{\alt}{\widetilde{\alpha}}
\newcommand{\nh}{\widehat{\nu}}
\newcommand{\Lh}{\widehat{L}}
\newcommand{\nt}{\widetilde{\nu}}
\newcommand{\Lt}{\widetilde{L}}
\newcommand{\Jt}{\widetilde{J}}
\newcommand{\lt}{\widetilde{\la}}
\newcommand{\Rt}{\widetilde{R}}
\newcommand{\rt}{\widetilde{r}}
\newcommand{\K}{K}
\newcommand{\Kn}{\sim_K}
\newcommand{\tK}{\widetilde{K}}
\newcommand{\Image}{\widetilde{Im}}

\section{Introduction}

Consider the Levi (block diagonal) subgroup $GL(\eta)\subset GL(n,\C)$
\begin{equation*}
  GL(\eta) = \prod_{i=1}^t GL(\eta_i,\C)
\end{equation*}
where $\eta=(\eta_1,\eta_2,\dots,\eta_t)$ is a sequence of positive
integers summing to $n$.  Define the
Littlewood-Richardson (LR) coefficient to be the multiplicity
\begin{equation}
	\LRC^\la_R = \dim \Hom_{GL(\eta)}(\Res^{GL(n)}_{GL(\eta)} V_\la,
	V_{R_1} \otimes V_{R_2} \otimes \dotsm \otimes V_{R_t})
\end{equation}
where $V_\la$ is the irreducible $GL(n,\C)$ module of highest weight $\la$
and $V_{R_i}$ is the irreducible $GL(\eta_i,\C)$-module of highest
weight $R_i$.  There is a well-known set $LRT(\la;R)$ of Young tableaux
(which shall be referred to as LR tableaux) whose cardinality is the
above coefficient $\LRC^\la_R$ \cite{LR}.

In \cite{SW} one of the authors and J. Weyman began the combinatorial
study of a family of polynomials $\K_{\la;R}(q)$ which are
$q$-analogues of the LR coefficients $\LRC^\la_R$ and are by definition
the Poincar\'e polynomials of isotypic components of Euler
characteristics of certain $\C[\gl_n]$-modules supported in nilpotent
conjugacy class closures.  The polynomials $\K_{\la;R}(q)$
were conjecturally described as the generating function
over \textit{catabolizable tableaux} with the charge statistic,
giving a simultaneous generalization of two
formulas of Lascoux and Sch\"utzenberger for the
Kostka-Foulkes polynomials \cite{La} \cite{LS1}.
The Kostka-Foulkes polynomials occur as special cases in two
different ways, namely, when each partition $R_i$ is a single row, or
each is a single column.  When $\mu$ has two columns, the 
Macdonald-Kostka polynomial $K_{\la,\mu}(q,t)$ has a nice
formula in terms of the polynomials $\K_{\la;R}(q)$ where each partition
$R_i$ has size at most $2$, \cite{Fi} \cite{St}.

Our point of departure is the observation that when each $R_i$ is
a \textit{rectangular} partition, the polynomial $\K_{\la;R}(q)$ seems
to coincide with another $q$-analogue of the appropriate LR coefficients,
given by the set $RC(\la;R)$ of \textit{rigged configurations} \cite{KR}.
One of the authors had already given a bijection
$\Psi_R:LRT(\la;R)\rightarrow RC(\la;R)$ \cite{K}.  The latter set is
endowed with a natural statistic $RC(\la;R)\rightarrow \NN$.  An obvious
problem is to give a direct description of the statistic on $LRT(\la;R)$
that is obtained by pulling back the statistic on $RC(\la;R)$ via the
bijection $\Psi_R$.  We offer two conjectures for this statistic,
which generalize the formulas for the charge statistic given by
Donin \cite{Do} and Lascoux, Leclerc, and Thibon \cite{LLT}.

Like the LR coefficients of which they are $q$-analogues,
the polynomials $\K_{\la;R}(q)$ satisfy symmetry and monotonicity
properties that extend those satisfied by the Kostka-Foulkes polynomials
\cite{SW}.  Indeed, some of these symmetries only appear after
generalizing from the Kostka-Foulkes case to the rectangular
LR case.  We give bijections and injections that exhibit these
properties combinatorially, for each of the three kinds of objects
(LR tableaux, catabolizable tableaux, and rigged configurations).
In particular, the monotonicity property is exhibited by
functorial statistic-preserving embeddings of families of LR tableaux,
generalizing the theory of the cyclage due to Lascoux and Sch\"utzenberger
\cite{LS3} \cite{La}.

There is another $q$-analogue of LR coefficients introduced
by Lascoux, Leclerc, and Thibon \cite{LLT} \cite{KLLT}.  Amazingly, these
polynomials arise in a completely different manner, namely, as
coefficient polynomials in a generating function over rim hook tableaux.  
We conjecture that the polynomials $\K_{\la;R}(q)$ coincide with theirs.

The paper is organized as follows.  Section \ref{gen func sec} recalls
the definition of the polynomial $\K_{\la;R}(q)$ and its
symmetry and monotonicity properties.  Sections \ref{RC section} through
\ref{CT section} give the three conjectured combinatorial descriptions
for the polynomials $\K_{\la;R}(q)$.  Section \ref{three bijections}
gives (conjecturally statistic-preserving) bijections
$\Psi_R$ from LR tableaux to rigged configurations
and $\Psi_{\rows(R)}$ from rigged configurations to
catabolizable tableaux.  For each of these kinds of objects
Sections \ref{symmetry section} through \ref{monotonicity section} give
bijections and injections that reflect the
symmetry and monotonicity properties of the polynomials $\K_{\la;R}(q)$.
These maps were defined so that they are intertwined by the
maps $\Psi_R$ and $\Psi_{\rows(R)}$.  

\section{Definition and properties of $\K_{\la;R}(q)$}
\label{gen func sec}
The material in this section essentially follows \cite{SW}.

\subsection{Generating function definition}
Let $\eta=(\eta_1,\eta_2,\dots,\eta_t)$ be a
sequence of positive integers that sum to $n$,
$\gamma=(\gamma_1,\gamma_2,\dots,\gamma_n)$ an integral weight
(that is, $\gamma\in\Z^n$), and
$\Roots_\eta$ the set of ordered pairs $(i,j)$ such that
$1\le i\le \eta_1+\eta_2+\dots+\eta_r < j\le n$ for some $r$.
Let $\Sigma_A$ be the symmetric group on the set $A$,
$[a,b]$ the closed interval of integers $i$ with $a\le i\le b$,
$[n]=[1,n]$, $x=(x_1,x_2,\dots,x_n)$ a sequence of variables,
$x^\gamma=x_1^{\gamma_1}x_2^{\gamma_2}\dots x_n^{\gamma_n}$,
and $\rho=(n-1,n-2,\dots,1,0)$.  The symmetric group $\Sigma_{[n]}$
acts on polynomials in $x$ by permuting variables.
Define the operators $J$ and $\pi$ by
\begin{equation*}
\label{J def}
\begin{split}
J(f) &= \sum_{w\in \Sigma_{[n]}} (-1)^w w(x^\rho f) \\
\pi f &= J(1)^{-1} J(f)
\end{split}
\end{equation*}
$J(1)$ is the Vandermonde determinant.
For the dominant (weakly decreasing) integral weight
$\la=(\la_1\ge \la_2\ge \dots\ge \la_n)$, the
character $s_\la(x)$ of $V_\la$ is given by the Laurent polynomial
\begin{equation*}
	s_\la(x) = \pi x^\la.
\end{equation*}
When $\la$ is a partition (that is, $\la_n\ge 0$), $s_\la$ is the
Schur polynomial.

Let $B_\eta(x;q)$, $H_{\gamma,\eta}(q)$, and
$\K_{\la,\gamma,\eta}(q)$ be the formal power series defined by
\begin{equation} \label{gen func}
\begin{split}
	B_\eta(x;q) &= \prod_{(i,j)\in \Roots_\eta}
		\dfrac{1}{1 - q x_i/x_j} \\
	H_{\gamma,\eta}(x;q) &= \pi (x^\gamma B_\eta(x;q)) \\
	&= \sum_\la s_\la(x) \K_{\la,\gamma,\eta}(q),
\end{split}
\end{equation}
where $\la$ runs over the dominant integral weights in $\Z^n$.
Ostensibly given by power series, the $\K_{\la,\gamma,\eta}(q)$ 
are in fact polynomials with integer coefficients \cite{SW}.

Let $R$ be the sequence $(R_1,R_2,\dots,R_t)$ with $R_i\in\Z^{\eta_i}$
a dominant integral weight for all $i$, and $\gamma(R)\in\Z^n$ 
the weight obtained by concatenating the parts of the $R_i$ in order.
Define
\begin{equation*}
  \K_{\la;R}(q) := \K_{\la,\gamma(R),\eta}(q).
\end{equation*}

It is known \cite{SW} that
\begin{equation} \label{LR spec}
	\K_{\la;R}(1) = \LRC^\la_R.
\end{equation}

From now on it is assumed that $\la$ is a partition and
each $R_i$ is a rectangular partition having $\eta_i$ rows
and $\mu_i$ columns.

\subsection{Special cases} Let $\la$ be a partition.
\begin{enumerate}
\item Let $R_i$ be the single row $(\mu_i)$ for all $i$,
where $\mu$ is a partition of length at most $n$.  Then
\begin{equation*}
  \K_{\la;R}(q) = K_{\la,\mu}(q),
\end{equation*}
the Kostka-Foulkes polynomial.  For a definition of the
Kostka-Foulkes polynomials (as well as the cocharge
and Macdonald-Kostka versions mentioned below) see \cite{Mac}.
\item Let $R_i$ be the single column $(1^{\eta_i})$ for all $i$.  Then
\begin{equation*}
\K_{\la;R}(q) = \tK_{\la^t,\eta^+}(q),
\end{equation*}
the cocharge Kostka-Foulkes polynomial, where $\la^t$ is the
\textit{transpose} or \textit{conjugate} of the partition $\la$
and $\eta^+$ is the partition obtained by sorting the parts of $\eta$
into weakly decreasing order.
\item Let $k$ be a positive integer and
$R_i$ the rectangle with $k$ columns and $\eta_i$ rows.
Then $\K_{\la;R}(q)$ is the Poincar\'e polynomial of the
isotypic component of the irreducible $GL(n)$-module
of highest weight $(\la_1-k,\la_2-k,\dots,\la_n-k)$ in
the coordinate ring of the Zariski closure
of the nilpotent conjugacy class whose Jordan form has
diagonal block sizes given by the transpose of the partition $\eta^+$.
In the case $\mu=(1^n)$ these are Kostant's generalized exponents
in type $A$.
\item Let $\la$ be a partition of $n$ and $\mu$
a two-column partition.  In this case J. Stembridge \cite{St} gave an
explicit formula for the Macdonald-Kostka polynomials, which has the form
\begin{equation*}
  \K_{\la,(2^r,1^{n-2r})}(q,t) =
  \sum_{k=0}^r q^k  \begin{bmatrix} r\\k \end{bmatrix}_t M^k_{r-k}(t)
\end{equation*}
where the $M^k_{r-k}(t)$ are members of a family of polynomials
$M^d_m(t)$ that are defined by iterated degree-shifted differences of
ordinary Kostka-Foulkes polynomials.  S. Fishel \cite{Fi} gave a
combinatorial description of the polynomials $M^d_m(t)$ in
terms of rigged configurations, using a variation of the original
statistic of \cite{KR} on the set of rigged configurations corresponding to 
standard tableaux.  Using the original statistic but replacing standard
tableaux by tableaux corresponding to sequences of tiny rectangles of the
form $(2)$, $(1,1)$, and $(1)$, we have
\begin{equation}\label{diff poly}
  M^k_{r-k}(t) = \K_{\la,((2)^{r-k},(1,1)^{k},(1)^{n-2r})} (t).
\end{equation}
The polynomials $M^d_m(t)$ are defined by a recurrence
that may be interpreted in terms of minimal degenerations
of nilpotent conjugacy class closures \cite{KKSW}. 
If the right hand side of \eqref{diff poly} is replaced by
a conjecturally equivalent formulation involving rigged configurations,
the resulting formula is proven in Section~\ref{monotonicity section}.
\end{enumerate}

\subsection{Defining recurrence}
The polynomials $\K_{\la;R}(q)$ satisfy the following recurrence
which generalizes Morris' recurrence for the Kostka-Foulkes polynomials
\cite{Mo} and Weyman's recurrence for the Poincar\'e polynomials in
case (3) above \cite{We}.  The polynomials $\K_{\la;R}(q)$ are uniquely
defined by the initial condition
\begin{equation*}
  \K_{\la;(R_1)}(q)=\delta_{\la,R_1}
\end{equation*}
and the recurrence
\begin{equation}
\label{recurrence}
  \K_{\la;R}(q) = 
  \sum_{w\in \Sigma_{[n]}/\Sigma_{[\eta_1]}\times\Sigma_{[\eta_1+1,n]}}
    (-1)^w q^{|\alpha(w)|}
      \sum_\tau \K_{\tau,\Rh}(q) \LRC^\tau_{\alpha(w),\beta(w)}
\end{equation}
where $w$ runs over the minimal length permutation in each coset in
$\Sigma_{[n]}$ of the given Young subgroup, 
$\alpha(w)$ and $\beta(w)$ are the first $\eta_1$ and last $n-\eta_1$ parts
of the weight $w^{-1}(\la+\rho)-(R_1+\rho)$, and $\Rh=(R_2,R_3,\dots,R_t)$.
The $w$-th summand is understood to be zero if $\alpha(w)$ has
a negative part.

The case that $R=(R_1,R_2)$ is now calculated explicitly.
Suppose $\mu_1\ge\mu_2$.  It is easy to show that for any $\la$,
\begin{equation} \label{mult free}
  \LRC^\la_{(R_1,R_2)} \in \{0,1\}.
\end{equation}
Assuming this multiplicity is one, it follows from \eqref{recurrence} that
\begin{equation} \label{two rectangles}
  \K_{\la;(R_1,R_2)}(q) = q^d
\end{equation}
where $d$ is the number of cells in $\la$ strictly to the
right of the $\mu_1$-th column.

\begin{ex} \label{P ex}
We give a running example.  Let $n=9$, $\mu=(3,2,1)$,
$\eta=(2,4,3)$.  Let 
\begin{equation*}
  R=((3,3),(2,2,2,2),(1,1,1))\qquad \la=(5,4,3,2,2,1,0,0,0)
\end{equation*}
Applying \eqref{recurrence}, all summands are zero except for
the identity permutation.  Using the LR rule, it is not hard to see that
the summand for $\tau$ is zero unless $\tau$ is one of the
three partitions $(3,3,3,2),$ $(3,3,2,2,1),$ and $(3,2,2,2,1,1)$.
Three applications of \eqref{two rectangles} yield
\begin{equation*}
\begin{split}
  \K_{\la;R}(q) &= q^3 (\K_{(3,3,3,2),\Rh}(q) +
  	2 \K_{(3,3,2,2,1),\Rh}(q) + \K_{(3,2,2,2,1,1),\Rh}(q)) \\
  	&= q^3 (q^3 + 2 q^2 + q^1) = q^6 + 2 q^5+ q^4.
\end{split}
\end{equation*}

\end{ex}

\subsection{Positivity}

Say that the sequence $R$ is \textit{dominant} if $\gamma(R)$ is.
The following positivity conjecture is due to Broer \cite{Broer}.

\begin{conj} \label{vanishing} \cite{Broer} If $R$ is dominant then
$\K_{\la;R}(q)\in\NN[q]$.
\end{conj}

This is known for some special cases, including the Kostka-Foulkes
case and the nilpotent adjoint orbit case (special cases (1) and
(3) in section \ref{gen func sec}).

\begin{ex} For $\la=(2,2)$ and $R=\{(1),(3)\}$, $\K_{\la;R}(q)=q-1$.
\end{ex}

\subsection{Cocharge normalization}

Given a sequence of rectangles $R$, let $r_{i,j}(R)$ be the number
of rectangles in $R$ that contain the cell $(i,j)$, or equivalently,
the number of indices $a$ such that $\mu_a \ge i$ and $\eta_a \ge j$.
Define the number
\begin{equation}
\label{n statistic}
  n(R) = \sum_{(i,j)} \binom{r_{i,j}(R)}{2}
\end{equation}
and the polynomial
\begin{equation}
\label{cocharge Poincare}
  \tK_{\la;R}(q) = q^{n(R)} \K_{\la;R}(q^{-1})
\end{equation}

When $R$ consists of single-rowed shapes, the above assertion
and both the two Kostka-Foulkes special cases imply that
\begin{equation*}
  K_{\la,\mu}(q) = q^{n(\mu)} \tK_{\la,\mu}(q^{-1})
\end{equation*}
where $n(\mu)=\sum_j \binom{\mu^t_j}{2}$.  This coincides with
the definition of the cocharge Kostka-Foulkes polynomials.

\begin{ex} \label{n example}
In our running example, the matrix $(r_{i,j})$ is given by
\begin{equation*}
(r_{i,j}) =
\begin{matrix}
3&2&1&0&\dots\\
3&2&1&0&\dots\\
2&1&0&0&\dots\\
1&1&0&0&\dots
\end{matrix}
\end{equation*}
with all other entries zero, so $n(R)=2 \binom{3}{2}+3 \binom{2}{2}=9$
and $\tK_{\la;R}(q) = q^5 + 2 q^4 + q^3$.
\end{ex}

Here is another version of the positivity conjecture,
which adds an observation about an upper bound on the
powers of $q$ that may occur.

\begin{conj}\label{con2} For $R$ dominant, $\tK_{\la;R}(q)\in\NN[q]$.
\end{conj}

\subsection{Reordering symmetry}

\begin{thm} \label{reordering} \cite{SW}
Suppose $R$ and $R'$ are dominant sequences of rectangles
that are rearrangements of each other.  Then
$\K_{\la;R}(q) = \K_{\la;R'}(q)$.
\end{thm}

This property is not obvious since the reordering of the sequence of
rectangles results in a major change in the generating function
$H_{\gamma,\eta}(x;q)$.

\subsection{Contragredient duality symmetry}
\label{duality}

Let $\rev(\eta)$ denote the reverse of $\eta$.
Fix a positive integer $k$ such that $k\ge\la_1$ and $k\ge\mu_i$
for all $i$.  Let $\lt=(k-\la_n,k-\la_{n-1},\dots,k-\la_1)$
and $\Rt_i=((k-\mu_i)^{\eta_i})$ for $1\le i\le t$.
Note that $\lt$ (resp. $\Rt_i$) is obtained by the 180 degree rotation of
the complement of the partition $\la$ (resp. $R_i$) inside the
$k \times n$ rectangle (resp. $k\times \eta_i$ rectangle).  Then

\begin{prop} \cite{SW} For $R$ dominant,
\begin{equation} \label{duality symmetry}
  \K_{\la;R}(q) = \K_{\lt;\rev(\Rt)}(q)
\end{equation}
\end{prop}

\begin{ex} In our running example, take $k=5$.  Then we have
\begin{equation*}
  \K_{(5,4,3,2,2,1,0,0,0),((3,3),(2,2,2,2),(1,1,1))}(q) =
  \K_{(5,5,5,4,3,3,2,1,0),((4,4,4),(3,3,3,3),(2,2))}(q).
\end{equation*}
\end{ex}

\subsection{Transpose symmetry}
Let $R^t$ be the sequence of rectangles
obtained by transposing each of the rectangles in $R$.

\begin{conj} \label{transpose symmetry}
Let $R$ be dominant and $R'$ a dominant rearrangement of $R^t$.  Then
\begin{equation}
\label{cocharge}
\K_{\la^t;R'}(q) = \tK_{\la;R}(q),
\end{equation}
where the left hand side is computed in $GL(m)$ where
$m$ is the total number of columns in the rectangles of $R$.
\end{conj}

This property is mysterious; it is not obvious from the
properties of the modules that define the polynomials
$\K_{\la;R}(q)$.

\begin{ex} In the running example, $\la^t=(6,5,3,2,1,0)$ and
\begin{equation*}
\begin{split}
  R^t&=((2,2,2),(4,4),(3)) \\
  R'&=((4,4),(3),(2,2,2))
\end{split}
\end{equation*}
Using \eqref{recurrence} we have
\begin{equation*}
  \K_{(6,5,3,2,1,0),((4,4),(3),(2,2,2))}(q)= q^5 + 2 q^4 + q^3,
\end{equation*}
which agrees with the computation of $\tK_{\la;R}(q)$ in
Example \ref{n example}.
\end{ex}

\subsection{Monotonicity}
\label{mono subsec}

Let $\alpha$ and $\beta$ be sequences of nonnegative integers.
Say that $\alpha\dom\beta$ if $\alpha^+ \dom\beta^+$ in the
dominance partial order on partitions.

Given a sequence of rectangles $R$, let $\tau^k(R)$
be the \textit{partition} whose parts consist of the $\eta_i$ 
such that $\mu_i=k$, sorted into decreasing order.
In other words, $\tau^k(R)$ is the multiset of heights of 
the rectangles in $R$ that have exactly $k$ columns.
Say that $R\dom R'$ if $\tau^k(R)\dom \tau^k(R')$ for all $k$.  

\begin{conj} \label{mono conj} Suppose $R$ and $R'$ are dominant
with $R\dom R'$.  Then
$\K_{\la;R'}(q) \ge \K_{\la;R}(q)$ coefficientwise.
\end{conj}

\begin{ex} Let $R$ be as usual and
$R'=((3),(3),(2,2,2,2),(1,1,1))$.  The sequence of partitions
$\tau^k(R)$ and $\tau^k(R')$ are given by
\begin{equation*}
\begin{split}
\tau^.(R) &= ((3),(4),(2),(),\dots) \\
\tau^.(R') &= ((3),(4),(1,1),(),\dots)
\end{split}
\end{equation*}
Now
\begin{equation*}
  \K_{\la;R'}(q) = 2 q^7 + 4 q^6 + 3 q^5 + q^4
\end{equation*}
which dominates $\K_{\la;R}(q)$ coefficientwise.
\end{ex}

Suppose all of the rectangles have the same number of columns $k$,
so that $\tau^j$ is empty for $j\not=k$.  In this case
Conjecture \ref{mono conj} may be verified as follows.
For the partition $\mu$ of $n$, let $A_\mu$ be the
coordinate ring of the closure of the nilpotent
conjugacy class of matrices with transpose Jordan type $\mu$.
If $\nu$ is another partition with $\mu\dom\nu$, then
restriction of functions gives a natural epimorphism
of graded $GL(n)$-modules $A_\nu\rightarrow A_\mu$.
This special case includes a monotonicity property of the
cocharge Kostka polynomials:
\begin{equation}
\label{Kostka mono}
  \tK_{\la,\nu}(q) \ge \tK_{\la,\mu}(q)\qquad \textit{if $\mu\dom\nu$.}
\end{equation}
The first combinatorial proof of this fact was given in \cite{La}
\cite{LS2} \cite {LS3}.  An alternate proof may be found in \cite{BK}.
It uses two substantial results.  The first is
the original interpretation \cite{LS1} of the Kostka-Foulkes
polynomial $\K_{\la,\mu}(q)$ as the generating function over
the set $CST(\la,\mu)$ of column-strict tableaux of shape
$\la$ and content $\mu$ with the charge statistic.  The second is the
statistic-preserving bijection in \cite{KR} that sends column-strict
tableaux to rigged configurations.  Under this map,
the image of $CST(\la,\mu)$ is a subset in the image of $CST(\la,\nu)$,
yielding \eqref{Kostka mono}.

 \subsection{Generalized Kostka polynomials and ribbon tableaux}
\label{gkp subsec}

According to \eqref{LR spec} and Conjecture~\ref{con2}, 
the generalized Kostka polynomials are $q$--analogues
of tensor product multiplicities.  Another $q$--analogue of 
tensor product multiplicities was introduced by A.~Lascoux, B.~Leclerc 
and J.-Y.~Thibon \cite{LLT} using the spin generating functions for the 
set of $p$--ribbon tableaux.  We refer the reader to
\cite[Sections~4 and 6]{LLT} for the definitions of the notions of 
a $p$--ribbon tableau $T$, the spin $s(T)$ of 
a $p$--ribbon tableau, and the ``$p$--ribbon version" of the modified 
Hall-Littlewood polynomials $\widetilde{G}_{\Lambda}^{(p)}(X_n;q)$.
Let $\Lambda$ be the partition with empty $p$--core and $p$-quotient
$(R_1,R_2,\dots,R_p)$ (see \cite[Chapter~1, Example~8]{Mac}, for the
definitions of the $p$--core and $p$--quotient of a partition $\lambda$).
By definition
\begin{equation*}
\widetilde{G}_{\Lambda}^{(p)}(X_n;q)=
\sum_{T\in\mathrm{Tab}_p(\Lambda,\le n)}q^{\widetilde{s}(T)}x^{wt(T)},
\end{equation*}
where the sum runs over the set $\mathrm{Tab}_p(\Lambda,\le n)$ of
$p$--ribbon tableaux of shape $\Lambda$ filled with numbers not exceeding
$n$, $wt(T)$ is the weight or content of the ribbon tableau $T$, and
$\widetilde{s}(T)=-s(T)+\max\{s(T)~|~T\in\mathrm{Tab}_p(\Lambda ,\le n)\}$
is the cospin of the $p$--ribbon tableau $T$, cf. \cite[(25)]{LLT}.
It is known \cite[Theorem~6.1]{LLT} that $\widetilde 
G_{\Lambda}^{(p)}(X_n;q)$ is a symmetric function. Following \cite{LLT},
let us define polynomials $\tK_{\Lambda\lambda}^{p)}(q)$ via the decomposition
\begin{equation*}
\widetilde{G}_{\Lambda}^{(p)}(X_n;q)=\sum_{\lambda}
	\tK_{\Lambda\lambda}^{(p)}(q) s_{\lambda}(X_n).
\end{equation*}
It is also well--known and goes back to D.~Littlewood \cite{StW}, that
$\tK_{\Lambda\lambda}^{(p)}(1)$ is the multiplicity of the highest 
weight $\lambda$ irreducible representation of $sl(n)$ in the tensor product 
$V_{R_1}\otimes \dots V_{R_p}$.

Recall that a sequence of partitions $R=(R_1,\ldots ,R_p)$ is 
called {\it dominant}, if for all $1\le i\le p-1$, the last part of
$R_i$ is at least as large as the first part of $R_{i+1}$.

\begin{conj}\label{kir1} Let $R=(R_1,\ldots ,R_p)$ be a dominant sequence of 
rectangular partitions, and $\Lambda$ the partition with 
empty $p$--core and $p$--quotient $(R_1,\ldots ,R_p)$. Then
\begin{equation*}
\tK_{\lambda;R}(q)=\tK_{\Lambda\lambda}^{(p)}(q).
\end{equation*}
\end{conj}

\section{Rigged configurations}
\label{RC section}

Rigged configurations are combinatorial objects that provide
a new interpretation of the Kostka-Foulkes
polynomial \cite{KR}.  However, this construction applies for
arbitrary sequences of rectangles, not just those consisting of all
single-rowed rectangles.  This section follows \cite{KR}.

Let $\la$ be a partition and $R$ a sequence of rectangular partitions
such that $|\la|=\sum_i |R_i|$.  Let $R_i$ have $\eta_i$ rows and
$\mu_i$ columns as usual.  A \textit{configuration} of type $(\la;R)$
is a sequence of partitions $\nu=(\nu^1,\nu^2,\dots)$ such that
\begin{equation*}
\begin{split} \label{config}
	|\nu^k| &= \sum_{j>k} \la_j - \sum_{a\ge1} \mu_a \max(\eta_a-k,0)\\
	&= - \sum_{j\le k} \la_j + \sum_{a\ge1} \mu_a \min(k,\eta_a)
\end{split}
\end{equation*}
for each $k\ge 1$.  Note that if $k\ge\ell(\la)$ and $k\ge\eta_a$ for all
$a$, then $\nu^k$ is empty.  We make the convention that
$\nu^0$ is the empty partition.

For a partition $\rho$, define the number
$Q_n(\rho) = \rho^t_1+\rho^t_2+\dots+\rho^t_n$, the number of cells
in the first $n$ columns of $\rho$.  The \textit{vacancy numbers}
of the configuration $\nu$ of type $(\la;R)$ are defined by
\begin{equation}\label{vacancies}
	P_{k,n}(\nu) = Q_n(\nu^{k-1})-2\, Q_n(\nu^k)+Q_n(\nu^{k+1})
		+ \sum_{a\ge1} \min(\mu_a,n)\, \delta_{\eta_a,k}
\end{equation}
for $k,n\ge 1$, where $\delta_{i,j}$ is the Kronecker delta.
Say that a configuration $\nu$ of type $(\la;R)$ is \textit{admissible}
if $P_{k,n}(\nu)\ge0$ for all $k,n\ge 1$.

\begin{ex} In the running example, there is a unique admissible
configuration of type $(\la;R)$, given by
\begin{equation*}
  \nu=((1),(2,1),(2,1),(2,1),(1))
\end{equation*}
Below is the table $P_{k,n}(\nu)$ of vacancy numbers for $k,n\ge 1$
with $k$ as row index and $n$ as column index.
\begin{equation*}
\begin{matrix}
0&1&1&\dots\\
0&0&1&\dots\\
1&1&1&\dots\\
0&0&0&\dots\\
0&1&1&\dots\\
1&1&1&\dots\\
0&0&0&\dots\\
\vdots&\vdots&\vdots&\dots
\end{matrix}
\end{equation*}
\end{ex}

A \textit{rigging} $L$ of an admissible configuration $\nu$ of type
$(\la;R)$ consists of an integer label $L^k_s$ for each row $\nu^k_s$ of
each of the partitions of $\nu$, such that:
\begin{enumerate}
\item $0 \le L^k_s \le P_{k,\nu^k_s}(\nu)$ for every $k\ge 1$ and
$1\le s\le \ell(\nu^k)$.
\item If $\nu^k_s = \nu^k_{s+1}$, then $L^k_s\ge L^k_{s+1}$.
\end{enumerate}
Only some of the vacancy numbers $P_{k,n}(\nu)$
appear as upper bounds for the labels $L^k_s$, namely, those where $n$ is
a part of the partition $\nu^k$.  The second condition merely says that
the labels for a fixed partition $\nu^k$ and part size $n$,
should be viewed as a multiset.

\begin{ex} To indicate the maximum rigging $\Lmax$ of the configuration
$\nu$, we replace each part $\nu^k_s$ by its maximum label
$P_{k,\nu^k_s}(\nu)$, obtaining
\begin{equation*}
  \Lmax=((0),(0,0),(1,1),(0,0),(0))
\end{equation*}
So $\nu$ has the four riggings
\begin{equation*}
\begin{split}
&((0),(0,0),(0,0),(0,0),(0))\qquad ((0),(0,0),(0,1),(0,0),(0)) \\
&((0),(0,0),(1,0),(0,0),(0))\qquad ((0),(0,0),(1,1),(0,0),(0))
\end{split}
\end{equation*}
\end{ex}

A \textit{rigged configuration} of type $(\la;R)$ is a pair $(\nu,L)$ where
$\nu$ is an admissible configuration of type $(\la;R)$ together with a
rigging $L$.  Let $C(\la;R)$ denote the set of admissible configurations
of type $(\la;R)$ and $RC(\la;R)$ the set of rigged configurations of type
$(\la;R)$.

The first important property of rigged configurations is that they
are enumerated by the LR coefficient $\LRC^\la_R$.

\begin{thm} \label{RC number} \cite{KR}
\begin{equation*}
	\LRC^\la_R = |RC(\la;R)| = \sum_{\nu\in C(\la;R)}
	\prod_{k,n\ge1} \binom{P_{k,n}(\nu)+m_n(\nu^k)}{m_n(\nu^k)}
\end{equation*}
where $m_n(\rho)$ denotes the number of parts of the partition $\rho$
of size $n$.
\end{thm}

Rigged configurations are endowed with a natural statistic we call
cocharge.  Let $(\nu,L)\in RC(\la;R)$.  Denote by
$\alpha_{k,n}=(\nu^k)^t_n$ 
the size of the $n$-th column of the $k$-th partition $\nu^k$ of $\nu$.
The \textit{cocharge} of $(\nu,L)$ is defined by
\begin{equation}
\label{RC cocharge}
\begin{split}
  \cocharge(\nu,L) &= \cocharge(\nu)+
     \sum_{k\ge 1}\sum_{s=1}^{\ell(\nu^k)} L^k_s \\
  \cocharge(\nu) &= \sum_{k,n\ge 1} \alpha_{k,n}(\alpha_{k,n}-\alpha_{k+1,n})
\end{split}
\end{equation}

Define the polynomial $RC_{\la;R}(q)$ by
\begin{equation*}
\begin{split}
  RC_{\la;R}(q) &= \sum_{(\nu,L)\in RC(\la;R)} q^{\cocharge{(\nu,L)}} \\
  &= \sum_{\nu\in C(\la;R)} q^{\cocharge(\nu)}
  	\prod_{k,n\ge1} \qbinom{P_{k,n}(\nu)+m_n(\nu^k)}{m_n(\nu^k)}
\end{split}
\end{equation*}
where $\qbinom{n}{k}$ is the $q$-binomial coefficient.

\begin{conj} For $R$ dominant,
\begin{equation}
	\tK_{\la;R}(q) = RC_{\la;R}(q)
\end{equation}
\end{conj}

\begin{ex} Let $\nu$ be the sole member of $RC(\la;R)$ as above.
We have
\begin{equation*}
\begin{split}
  \cocharge(\nu) &= 1(1-2)+2(2-2)+2(2-2)+2(2-1)+1(1-0) + \\
  	& 0(0-1)+1(1-1)+1(1-1)+1(1-0) \\
  	&= (-1+2+1)+(-1+1) = 3,
\end{split}
\end{equation*}
so that $RC_{\la;R}(q) = q^3(1+q)(1+q)$,
which agrees with $\tK_{\la;R}(q)$.

\end{ex}

\begin{comment}
\begin{ex} There is a single admissible configuration of type $(\la^t;R^t)$
given by
\begin{equation*}
  \nh=((3),(3,1),(2,1),(1)).
\end{equation*}
The transposes of these partitions
are $((1,1,1),(2,1,1),(2,1),(1))$.  We have
\begin{equation*}
\begin{split}
  \cocharge(\nh)&= (1(1-2)+2(2-2)+2(2-1)+1(1-0))\\
  	&+(1(1-1)+1(1-1)+1(1-0)) + (1(1-1)+1(1-0)) \\
  &= (-1+0+2+1)+(0+0+1)+(0+1) = 4.
\end{split}
\end{equation*}
The table of vacancy numbers $P_{k,n}(\nh)$ is given by
\begin{equation*}
\begin{matrix}
1&1&1&1&\dots\\
0&1&1&2&\dots\\
0&0&1&1&\dots\\
0&1&1&1&\dots\\
1&1&1&1&\dots\\
0&0&0&0&\dots\\
\vdots&\vdots&\vdots&\vdots&\dots
\end{matrix}
\end{equation*}
The maximum rigging for $\nh$ is $((1),(1,0),(0,0),(0))$.  There are
four riggings for $\nh$, namely,
\begin{equation*}
\begin{split}
&((1),(1,0),(0,0),(0))\qquad((1),(0,0),(0,0),(0)) \\
&((0),(1,0),(0,0),(0))\qquad((0),(0,0),(0,0),(0))
\end{split}
\end{equation*}
Thus $RC_{\la^t;R^t}(q) = q^6 + 2 q^5+q^4$, which agrees with
Example \ref{P ex}.
\end{ex}
\end{comment}

\section{Littlewood-Richardson tableaux}
\label{LR section}

This section discusses a set $LRT(\la;R)$ of tableaux that has cardinality
$\LRC^\la_R$.

Let us fix notation.  Let the intervals $A_1$, $A_2$, etc., be given by
dividing the interval $[n]$ into successive subintervals
of sizes $\eta_1$, $\eta_2$, etc.  Let $K_i$ be the rectangular
tableau of shape $R_i$ whose $j$-th row consists of $\mu_i$ copies of the
$j$-th largest letter of the interval $A_i$ for $1\le j\le \eta_i$.

\begin{ex} \label{notation ex}
Let $\lambda=(5,4,3,2,2,1)$ and $R=(R_1,R_2,R_3)$ with
$R_1=(3,3)$, $R_2=(2,2,2,2)$ and $R_3=(1,1,1)$.
$\mu=(3,2,1)$, $\eta=(2,4,3)$,
$\gamma=(3,3,2,2,2,2,1,1,1)$,
$A_1=[1,2]$, $A_2=[3,6]$, $A_3=[7,9]$.
\begin{equation*}
K_1 = \begin{matrix} 1&1&1\\2&2&2 \end{matrix} \qquad
K_2 = \begin{matrix} 3&3\\4&4\\5&5\\6&6 \end{matrix} \qquad
K_3 = \begin{matrix} 7\\8\\9 \end{matrix}
\end{equation*}
\end{ex}

Say that a word is \textit{lattice} if the content of every final
subword is a partition.  A (skew) tableau is said to be lattice if its
row-reading word is lattice (see the example below).

Let $LRT(\la;R)$ be the set of column
strict tableaux $T$ of shape $\lambda$ and content $\gamma=\gamma(R)$
such that for all $i$, the restriction $T|_{A_i}$ of the tableau $T$ to the
alphabet $A_i$ is lattice with respect to the alphabet $A_i$, or
equivalently, is Knuth equivalent to the tableau $K_i$.
We call these LR tableaux.
Note that if each rectangle of $R$ is a single row, then
$LRT(\la;R)=CST(\la,\gamma)$.

\begin{ex} \label{LR example}
$LRT(\la;R)$ consists of the following four tableaux.
\begin{equation}
\begin{matrix}
1&1&1&3&3\\
2&2&2&4& \\
4&5&7& & \\
5&6& & & \\
6&8& & & \\
9& & & &
\end{matrix}\qquad
\begin{matrix}
1&1&1&3&3\\
2&2&2&7& \\
4&4&8& & \\
5&5& & & \\
6&6& & & \\
9& & & &
\end{matrix}\qquad
\begin{matrix}
1&1&1&3&7\\
2&2&2&4& \\
3&5&8& & \\
4&6& & & \\
5&9& & & \\
6& & & &
\end{matrix}\qquad
\begin{matrix}
1&1&1&3&7\\
2&2&2&8& \\
3&4&9& & \\
4&5& & & \\
5&6& & & \\
6& & & &
\end{matrix}
\end{equation}
Let $T$ be the first tableau.  The skew tableaux $T|_{A_i}$ are given below.
\begin{equation*}
T|_{A_1}=\begin{matrix}
1&1&1&\cdot&\cdot\\
2&2&2&\cdot& \\
\cdot&\cdot&\cdot& & \\
\cdot&\cdot& & & \\
\cdot&\cdot& & & \\
\cdot& & & &
\end{matrix}\qquad
T|_{A_2}=\begin{matrix}
\cdot&\cdot&\cdot&3&3\\
\cdot&\cdot&\cdot&4& \\
4&5&\cdot& & \\
5&6& & & \\
6&\cdot& & & \\
\cdot& & & &
\end{matrix}\qquad
T|_{A_3}=\begin{matrix}
\cdot&\cdot&\cdot&\cdot&\cdot\\
\cdot&\cdot&\cdot&\cdot& \\
\cdot&\cdot&7& & \\
\cdot&\cdot& & & \\
\cdot&8& & & \\
9& & & &
\end{matrix}
\end{equation*}
We have $\word(T|_{A_1})=222111$, $\word(T|_{A_2})=65645433$, and
$\word(T|_{A_3})=987$, which are lattice in their
respective alphabets.
\end{ex}

The following result is an easy consequence of the classical
Littlewood-Richardson rule \cite{LR}.

\begin{prop} $|LRT(\la;R)| = \LRC^\la_R$.
\end{prop}

Next we define a statistic $\charge_R$ on LR tableaux.

Say that a word $v$ is \textit{$R$-LR} (short for $R$-Littlewood-Richardson)
if $P(v)\in LRT(\la;R)$
for some $\la$, or equivalently, that $v|_{A_i}$ is lattice in the alphabet
$A_i$ for all $i$, or equivalently, that $P(v|_{A_i})=K_i$ for all $i$.
Clearly an $R$-LR word must have content $\gamma(R)$.

Let $v$ be a word of partition content.
Say that the sequence of words $\{v^1,v^2,\dots\}$ is a
\textit{standard decomposition} of $v$, if the $v^i$ are
standard words of weakly decreasing size and
and $v$ is a shuffle of the words $\{v^i\}$, that is,
the $v^i$ may be simultaneously disjointly embedded in $v$
and exhaust the letters of $v$.

Now let $R$ be dominant, $v$ an $R$-LR word,
and $\{v^i\}$ a standard decomposition of $v$.  Say that
$\{v^i\}$ is \textit{proper} if $v^i|_{A_j}$ is either
empty or the decreasing word consisting of the letters of $A_j$,
for all $i$ and $j$.  Assuming this holds,
let $u^i$ be the reverse of the word obtained from $v^i$ by
replacing each letter in $A_j$ by the letter $j$.
Define
\begin{equation*}
  \charge_R(v) = \min_{\{v^i\}} \sum_i \cocharge(u^i)
\end{equation*}
where $\{v^i\}$ runs over the proper standard decompositions of $v$
and $\{u^i\}$ is related to $v^i$ as above.
Say that a proper standard decomposition of an $R$-LR word $v$
is minimal if it attains the minimum in the definition of
$\charge_R(v)$.

\begin{enumerate}
\item When $\mu$ is a partition and $R_i=(\mu_i)$ for all $i$,
the subalphabet $A_j=\{j\}$, so every standard decomposition
is proper and $u^i$ is the reverse of $v^i$.
But for standard words, the cocharge of the reverse of a word
equals the charge.  So $\charge_R$ is the usual charge statistic,
since Donin \cite{Do} asserts that the particular standard
decomposition given by Lascoux and Sch\"utzenberger \cite{LS1}
is always minimal.
\item When $R_i=(1^{\eta_i})$ for all $i$,
an $R$-LR word $v$ is necessarily standard, so $v^1=v$
and by definition $\charge_R(v)=\cocharge(u^1)$.
\end{enumerate}

\begin{ex} Let
$R=((2,2),(2,2,2),(2),(1))$ so that $\mu=(2,2,2,1)$ and $\eta=(2,3,1,1)$.
The alphabets $A_j$ are given by $[1,2]$, $[3,5]$, $[6,6]$ and $[7,7]$.
Let $\la=(5,4,2,1,1,0,0)$.  We indicate part of the computation of
$\charge_R$ for the tableau $T$ given below.
\begin{equation*}
T=\begin{matrix}
1&1&3&3&6\\
2&2&4&6& \\
4&5& & & \\
5& & & & \\
7& & & &
\end{matrix}
\end{equation*}
We have $\word(T)=7.5.45.2246.11336$.  It turns out that
$\charge_R(\word(T))=7$.  There are at most
$2^3=8$ distinct proper standard decompositions of $\word(T)$ depending on
choices involving the letters $4$, $5$, and $6$; the pairs of ones,
twos, and threes
occur side by side and thus do not generate other proper decompositions.
The unique minimal decomposition is $v^1=7542613$, $v^2=524136$.
To see this, reversing these subwords gives $(3162457,631425)$.
Replacing each letter in $A_j$ by $j$ we have the pair of words
$(u^1,u^2)=(2131224,321212)$.  To take the cocharge we act by automorphisms
of conjugation to change them to partition content, obtaining
$(2131124,321211)$.  Taking the cocharges of these words,
we obtain $(3,4)$ whose sum is $7$.  Let us now perform the same computation
for the circular standard decomposition of Lascoux and
Sch\"utzenberger (which is always proper): $v^1=7524136$ and $v^2=542613$.
We have $u^1=3212124$ and $u^2=213122$.  Acting by automorphisms of
conjugation to move to partition content, we have $(3212114,213112)$, whose
cocharges are $(6,2)$ whose sum is $8$, which is not minimal.
\end{ex}

It would be desirable to have an algorithm that computes
a minimal standard decomposition.  The above example shows that
the algorithm of Lascoux and Sch\"utzenberger for selecting
a standard decomposition in the computation of charge,
does not work for $\charge_R$.

For $R$ dominant, define
\begin{equation*}
  LRT_{\la;R}(q) = \sum_{T\in LRT(\la;R)} q^{\charge_R(T)}.
\end{equation*}

\begin{conj} \label{K as LRT} For $R$ dominant,
\begin{equation*}
  \K_{\la;R}(q) = LRT_{\la;R}(q).
\end{equation*}
\end{conj}

\section{Catabolizable tableaux}
\label{CT section}

We recall from \cite{SW} the notion of an $R$-catabolizable tableau.
Recall the subintervals $A_i$ and the canonical rectangular tableaux $K_i$
from the definition of $LRT(\la;R)$.  Let $S$ be a column strict tableau of
partition shape.  Suppose that $S|_{A_1}=K_1$.
Let $S_+$ be the first $\eta_1$ rows of the skew tableau $S-K_1$
and let $S_-$ be the remainder.
Define the (row) \textit{$R_1$-catabolism} $\cat_{R_1}(S)$ of $S$
to be the tableau $P(S_+ S_-)$, where $P(w)$ is the Schensted insertion
tableau for the word $w$ \cite{Sch}.  Say that $S$ is
\textit{$R$-catabolizable} if $S|_{A_1}=K_1$ and $\cat_{R_1}(S)$ is
$\Rh$-catabolizable in the alphabet $[\eta_1+1,n]=[n]-A_1$,
where $\Rh=(R_2,R_3,\dotsc)$.  The empty tableau is
considered to be the unique catabolizable tableau for the empty sequence
of rectangles.  Note that an $R$-catabolizable tableau must have content
$\gamma(R)$.

Denote by $CT(\la;R)$ the set of $R$-catabolizable tableaux of shape $\la$.

\begin{ex} \label{cat ex}
Recall the tableaux $K_1$, $K_2$, and $K_3$ from Example
\ref{notation ex}.  Here are the four tableaux that comprise the
set $CT(\la;R)$.
\begin{equation*}
\begin{matrix}
1&1&1&5&6\\
2&2&2&6& \\
3&3&7& & \\
4&4& & & \\
5&8& & & \\
9& & & &
\end{matrix}\qquad
\begin{matrix}
1&1&1&6&6\\
2&2&2&7& \\
3&3&8& & \\
4&4& & & \\
5&5& & & \\
9& & & &
\end{matrix}\qquad
\begin{matrix}
1&1&1&5&7\\
2&2&2&6& \\
3&3&8& & \\
4&4& & & \\
5&9& & & \\
6& & & &
\end{matrix}\qquad
\begin{matrix}
1&1&1&6&7\\
2&2&2&8& \\
3&3&9& & \\
4&4& & & \\
5&5& & & \\
6& & & &
\end{matrix}
\end{equation*}
Let $S$ be the first of these tableaux.  The following calculation
shows that $S$ is $R$-catabolizable.
\begin{equation*}
\begin{split}
S&=\begin{matrix}
1&1&1&5&6\\
2&2&2&6& \\
3&3&7& & \\
4&4& & & \\
5&8& & & \\
9& & & &
\end{matrix}\qquad
S-K_1 = \begin{matrix}
\cdot&\cdot&\cdot&5&6\\
\cdot&\cdot&\cdot&6& \\
3&3&7& & \\
4&4& & & \\
5&8& & & \\
9& & & &
\end{matrix} \\
S_+ S_- &=
\begin{matrix}
 & &3&3&7\\
 & &4&4& \\
 & &5&8& \\
 & &9& & \\
5&6& & & \\
6& & & &
\end{matrix} \\
\cat_{R_1}(S)&=
\begin{matrix}
3&3&7\\
4&4&8\\
5&5&9\\
6&6&
\end{matrix} \qquad
\cat_{R_2} \cat_{R_1} S =
\begin{matrix}
7\\8\\9
\end{matrix}
\qquad
\cat_{R_3} \cat_{R_2} \cat_{R_1} S = \varnothing
\end{split} 
\end{equation*}
\end{ex}

Define 
\begin{equation*}
  CT_{\la;R}(q) = \sum_{S\in CT(\la;R)} q^{\charge(S)}.
\end{equation*}

\begin{conj} \label{P as CT} Let $R$ be dominant.  Then
\begin{equation*}
\K_{\la;R}(q) = CT_{\la;R}(q)
\end{equation*}
\end{conj}

When each rectangle in $R$ is a single row, this formula
(as well as that in Conjecture \ref{K as LRT}) specializes
to the description of the Kostka-Foulkes polynomials in \cite{LS1}.

When each rectangle is a single column, this formula specializes
(with some work \cite{SW}) to a formula for the cocharge
Kostka-Foulkes polynomials \cite{La}.

\section{Bijections between the three sets}
\label{three bijections}

In \cite{KR} two bijections $LRT(\la;R)\rightarrow RC(\la;R)$ were
given, in the case that all were single rows or all single columns
When all are single rows, $LRT(\la;R)=CST(\la,\mu)$,
the set of column strict tableaux of shape $\la$ and content $\mu$.
When all rectangles are single columns, there is an obvious bijection
$LRT(\la;R)\rightarrow CST(\la^t,\eta)$ given by taking the transpose
of a tableau and then relabelling (see the LR-transpose map in
Section \ref{transpose section}).  In \cite{KR} it is asserted that
the first of these two bijections is statistic-preserving.
One of the authors \cite{K} has given a common generalization of these
bijections, which we denote by
\begin{equation*}
\Phi_R:LRT(\la;R)\rightarrow RC(\la^t;R^t).
\end{equation*}
For those that have some familiarity with \cite{KR}
we point out two twists in its definition.
The first difference is in the labelling convention.
In \cite{KR} the bijections use the notion of a \textit{singular
string} in a rigged configuration $(\nu,J)$, that is, a part $\nu^k_s$
whose label $J^k_s$ attains the maximum value $P_{k,\nu^k_s}(\nu)$.
This convention is called the quantum number labelling.
Here we employ the coquantum number labelling, in which a 
singular string is a part $\nu^k_s$ whose label $L^k_s$ is zero.
The second difference is the direction in which the
cells of the rectangles $R_i$ are ordered.  In \cite{KR}
the cells of the one-rowed rectangles are ordered along rows,
but here the cells of the rectangles are ordered along columns.
This accounts for the transposing of shapes in passing from
LR tableaux to rigged configurations.

\begin{conj} \label{LR to RC preserving} For $R$ dominant,
the bijection $\Phi_R$ is statistic-preserving, that is,
for $T\in LRT(\la;R)$ and $(\nu,L)=\Phi_R(T)$, we have
\begin{equation*}
  \charge_R(T) = \cocharge(\nu,L)
\end{equation*}
\end{conj}

The dominance of $R$ is only necessary for the definition of $\charge_R$.

This bijection may also be used to define a map from rigged configurations
to catabolizable tableaux.  Let $\rows(R)$ be the sequence of 
rectangles obtained by slicing each rectangle of $R$ into single rows.
Recall the weight $\gamma=\gamma(R)$; its parts give the single-rowed
shapes of $\rows(R)$.  In our running example,
$\rows(R)=((3),(3),(2),(2),(2),(2),(1),(1),(1))$
and $\gamma=(3,3,2,2,2,2,1,1,1)$.
Similarly, define $\cols(R^t)$ to be the sequence of transposes of
$\rows(R)$.  We have the bijection
\begin{equation*}
  \Phi_{\rows(R)}^{-1}:
RC(\la^t;\cols(R^t))\rightarrow LRT(\la,\rows(R))=CST(\la,\gamma)
\end{equation*}
It is clear from the definitions that there is an
inclusion 
\begin{equation*}
  RC(\la^t;R^t)\subseteq RC(\la^t;\cols(R^t)).
\end{equation*}
Let
\begin{equation*}
b_R:RC(\la^t;R^t)\rightarrow CST(\la,\gamma)
\end{equation*} 
be the restriction of the map $\Phi_{\rows(R)}^{-1}$ to the subset
$RC(\la^t;R^t)$.

\begin{comment}
Say that an $a\times b$ rectangle \textit{swallows}
a $c\times d$ rectangle if $a>c$ and $b>d$.  Say that the sequence
of rectangles \textit{respects swallowing} if no rectangle
swallows a previous one.
\end{comment}

\begin{conj} \label{image conj}
\begin{comment}
Let $R$ be a sequence of rectangles that respects swallowing.
\end{comment}
If $R$ is dominant then $\Image b_R = CT(\la;R)$, and $b_R \circ \Phi_R$
is a bijection $LRT(\la;R)\rightarrow CT(\la;R)$
sending $\charge_R$ to $\charge$.
\end{conj}

In section \ref{monotonicity section} the composite
map $LRT(\la;R)\rightarrow CT(\la;R)$ is given without
using rigged configurations as an intermediate set.

\section{Symmetry Bijections}
\label{symmetry section}

Let $R$ and $R'$ be any two sequences of rectangles that are
rearrangements of each other.  It follows immediately from the definitions
that $RC(\la;R)=RC(\la;R')$, so that
\begin{equation*}
  RC_{\la;R}(q)=RC_{\la;R'}(q)
\end{equation*}

We now give the symmetry bijections of LR tableaux and a
another conjectural description for $\charge_R$.
Fix the sequence of rectangles $R=(R_1,R_2,\dots,R_t)$.
Let $\SR$ be the set of rearrangements of $R$.  $\Sigma_{[t]}$ acts on
$\SR$ in the obvious way.  For $u\in \Sigma_{[t]}$,
we wish to define bijections
\begin{equation*}
  u_R:LRT(\la;R)\rightarrow LRT(\la;u R)
\end{equation*}
that give an action of $\Sigma_{[t]}$ on LR tableaux
in the sense that the two bijections
$LRT(\la;R)\rightarrow LRT(\la;v u R)$ given by
$(v \circ u)_R$ and $v_{u R} \circ u_R$,
coincide.  When each rectangle in $R$ is a single row,
these bijections coincide with the action of the symmetric group on the
plactic algebra by the automorphisms of conjugation \cite{LS3},
whose Coxeter generators are sometimes called crystal reflection operators.

Suppose that $u$ is the adjacent transposition $s_1=(1 2)$
and $R=(R_1,R_2)$ consists of two rectangles.  Let $P\in LRT(\la;R)$.
By \eqref{mult free},
both the sets $LRT(\la;R)$ and $LRT(\la;(R_2,R_1))$ are singletons,
say $\{P\}$ and $\{P'\}$ respectively.
In this case the bijection $s_1=(s_1)_R$ is defined by $s_1 P=P'$, and its
inverse (also denoted $s_1$ by suppressing the subscript) is given by
$s_1 P'=P$.  The tableau $P'$ can be calculated as follows.
Let $A_1'$, $A_2'$, $K_1'$ and $K_2'$ be the subalphabets and canonical
tableaux for the sequence of rectangles $s_1 R=(R_2,R_1)$.
Clearly $P'|_{A_1'}=K_1'$.  The remainder $P'|_{A_2'}$ of $P'$ is the
column strict tableau of shape $\la/R_2$ obtained by the jeu-de-taquin
given by sliding the tableau $K_2'$ to the southeast into the skew
shape $\la/R_1$ using the order of cells defined by the skew tableau
$P|_{A_2}$.

Next suppose $u=s_p$ where $p>1$.  Let $B=A_p\cup A_{p+1}$.
Let $(P,Q)$ be the tableau pair corresponding to the column
insertion of the row-reading word of the skew tableau
$T|_B$.  Since $T\in LRT(\la;R)$, the tableaux
$T|_{A_p}$ and $T|_{A_{p+1}}$ are lattice with
respect to the alphabets $A_p$ and $A_{p+1}$ respectively.
Since the lattice condition is invariant under Knuth equivalence,
it follows that the restrictions of $P$ to the alphabets
$A_p$ and $A_{p+1}$ are lattice, so that $P\in LRT(\rho;(R_p,R_{p+1}))$
in the alphabet $B$, where $\rho$ is the shape of $P$.
Let $P'$ be the unique tableau in the singleton set
$LRT(\rho;(R_{p+1},R_p))$ in the alphabet $B$.  By \cite{Wh}[Theorem 1],
there is a column strict tableau $U$ of the same skew shape as $T|_B$,
whose row-reading word corresponds to the tableau pair $(P',Q)$
under column insertion.  Define $s_p T$ be the tableau which agrees with
$U$ on the alphabet $B$ and agrees with $T$ on the complement of $B$.
Note that $s_p T\in LRT(\la;s_p R)$ since latticeness is preserved
under Knuth equivalence.

\begin{ex} \label{exchange rectangles}
Let $p=2$.  The subalphabets for the sequence of rectangles
$s_p R$ are $A_1'=[1,2]$, $A_2'=[3,5]$, and $A_3'=[6,9]$, with
$B=[3,9]$.  Consider the first tableau $T$ of Example \ref{LR example}.
The tableau $s_2 T$ is computed as follows.
\begin{equation*}
\begin{split}
T&=\begin{matrix}
1&1&1&3&3\\
2&2&2&4& \\
4&5&7& & \\
5&6& & & \\
6&8& & & \\
9& & & &
\end{matrix} \qquad
T|_B=\begin{matrix}
\cdot&\cdot&\cdot&3&3\\
\cdot&\cdot&\cdot&4& \\
4&5&7& & \\
5&6& & & \\
6&8& & & \\
9& & & &
\end{matrix} \\
\word(T) &= 9.68.56.457.4.33\\
P&=\begin{matrix}
3&3&7\\
4&4& \\
5&5& \\
6&6& \\
8& & \\
9& & 
\end{matrix} \qquad
Q=\begin{matrix}
1&2&6\\
3&5& \\
4&8& \\
7&10& \\
9& & \\
11& &
\end{matrix} \qquad
P'=\begin{matrix}
3&6&6\\
4&7& \\
5&8& \\
7&9& \\
8& & \\
9& & 
\end{matrix} \\
(s_2 T)|_B&=\begin{matrix}
\cdot&\cdot&\cdot&3&6\\
\cdot&\cdot&\cdot&4& \\
5&6&7& & \\
7&8& & & \\
8&9& & & \\
9& & & &
\end{matrix} \qquad
s_2 T=\begin{matrix}
1&1&1&3&6\\
2&2&2&4& \\
5&6&7& & \\
7&8& & & \\
8&9& & & \\
9& & & &
\end{matrix}
\end{split}
\end{equation*}
\end{ex}

For an arbitrary permutation $\sigma\in \Sigma_{[t]}$, the bijection
\begin{equation*}
u_R:LRT(\lambda;R)\rightarrow LRT(\lambda;u R)
\end{equation*}
is defined to be the composition $u=s_{a_1} s_{a_2} \dots s_{a_r}$,
where $a_1 a_2 \dots a_r$ is a reduced word for $u$.
These bijections were chosen with the following property in mind.

\begin{conj} \label{LRT symmetry} The following diagram commutes:
\begin{equation*}
\begin{CD}
	LRT(\la;R) @>{u_R}>> LRT(\la;u R) \\
	@V{\Phi_R}VV		@VV{\Phi_{u R}}V	\\
	RC(\la^t;R^t) @= RC(\la^t;u (R^t)) 
\end{CD}
\end{equation*}
In particular, the bijection $u_R$ is independent of the reduced word
of $u$, and the bijections of the form $u_R$ define an action of
$\Sigma_{[t]}$ on the collection $\SR$ of LR tableaux.
\end{conj}

When each rectangle in $R$ is a single row, 
a version of this result is stated in \cite{BK}[(2.17)].

It suffices to prove Conjecture \ref{LRT symmetry} for adjacent
transpositions $s_p$.  We show that this conjecture may be reduced
to another conjecture on evacuation.
Let $\ev=\ev_{[n]}$ denote the evacuation involution on column strict
tableaux in the alphabet $[n]$; it is defined by the conditions
\begin{equation*}
  \shape(\ev(T)|_{[k]}) = \shape(P(T|_{[n+1-k,n]}))
\end{equation*}
for all $1\le k\le n$.  Since latticeness is preserved by Knuth equivalence,
$\ev$ restricts to a bijection
$LRT(\la;R) \rightarrow LRT(\la;\rev(R))$ where $\rev(R)$ is the reverse
of $R$.

Consider also the involution $\theta$ on
$RC(\la;R)=RC(\la;\rev(R))$ given by $(\nu,L)\mapsto (\nu,J)$ where
\begin{equation*}
  L^k_s + J^k_s = P_{k,\nu^k_s}(\nu)
\end{equation*}
for all $k\ge1$, $1\le s\le \ell(\nu^k)$.  This is slightly sloppy since
the labels have to be reordered to satisfy the formal definition
of a rigging.  The involution $\theta$ complements each coquantum
number $L^k_s$ with respect to its maximum possible value.

\begin{conj} \label{ev conj} The following diagram commutes:
\begin{equation*}
\begin{CD}
	LR(\la;R) @>\ev>> LR(\la;\rev(R)) \\
	@V{\Phi_R}VV		@VV{\Phi_{\rev(R)}}V	\\
	RC(\la^t;R^t) @>>\theta> RC(\la^t;\rev(R^t)) 
\end{CD}
\end{equation*}
\end{conj}

When each rectangle of $R$ is a single row, a similar assertion
is made in \cite{KR}.

\begin{lem} Conjecture \ref{LRT symmetry} follows from
Conjecture \ref{ev conj}.
\end{lem}
\begin{proof} We may assume that $u=s_p$.  If $p=1$ this may be verified
directly using the definition of the map $\Phi_R$.  So suppose $p>1$.
Again by the definition of $\Phi_R$, 
it suffices to assume that $t=p+1$, that is, $s_p$ exchanges
the last two rectangles in $R$.  We have
\begin{equation*}
  \ev s_p T = s_1 \ev T
\end{equation*}
for all $T\in LRT(\la;R)$, which holds since
latticeness is preserved by Knuth equivalence.  Using the
fact that $\ev$ is an involution, that Conjecture \ref{LRT symmetry}
holds for $s_1$, and assuming Conjecture \ref{ev conj}, it follows
that Conjecture \ref{LRT symmetry} also holds for $s_p$.
\end{proof}

Next is another conjectural characterization of the statistic
$\charge_R$ on LR tableaux; its definition requires  
Conjecture \ref{LRT symmetry}.  The notation in the definition of the
bijection $s_p$ is used here.  Let $T\in LRT(\la;R)$.  Define the statistic 
$d_{p,R}(T)$ to be the number of cells in $P=P(T|_{A_p\cup A_{p+1}})$
that lie to the right of the $c$-th column, where
$c=\max(\mu_p,\mu_{p+1})$.  Compare this to \eqref{two rectangles}.  Let
\begin{equation*}
  d_R(T) = \sum_{p=1}^{t-1} (t-p)\, d_{p,R}(T).
\end{equation*}
The alternate definition for $\charge_R$ (valid for any $R$) is:
\begin{equation} \label{charge_R}
  \charge_R(T) = 1/t! \sum_{u \in\Sigma_{[t]}} d_{u R}(u T)
\end{equation}
(see Conjecture \ref{K as LRT}).  It is not clear why this
quantity should be an integer.  When each rectangle in $R$ is a single row,
\eqref{charge_R} specializes to the formula for the charge given in
\cite{LLT}.  By definition 
\begin{equation}
  \charge_{u R}(u T) =\charge_R(T)
\end{equation}
for all $T\in LRT(\la;R)$ and $u\in\Sigma_{[t]}$.

\begin{ex} Let us use \eqref{charge_R} to
compute $\charge_R$ for the tableau $T$ of Example
\ref{exchange rectangles}.  It is necessary to compute the entire
orbit of $T$ under the action of $\Sigma_{[3]}$.
\begin{equation*}
\begin{split}
T&=\begin{matrix}
1&1&1&3&3\\
2&2&2&4& \\
4&5&7& & \\
5&6& & & \\
6&8& & & \\
9& & & &
\end{matrix} \\
s_1 T&=\begin{matrix}
1&1&5&5&5\\
2&2&6&6& \\
3&3&7& & \\
4&4& & & \\
6&8& & & \\
9& & & &
\end{matrix} \qquad
s_2 T =\begin{matrix}
1&1&1&3&6\\
2&2&2&4& \\
5&6&7& & \\
7&8& & & \\
8&9& & & \\
9& & & &
\end{matrix} \\
s_2 s_1 T &= \begin{matrix}
1&1&5&8&8\\
2&2&6&9& \\
3&3&7& & \\
4&4& & & \\
8&9& & & \\
9& & & &
\end{matrix} \qquad
s_1 s_2 T = \begin{matrix}
1&4&4&4&6\\
2&5&5&5& \\
3&6&7& & \\
7&8& & & \\
8&9& & & \\
9& & & & 
\end{matrix} \\
s_1 s_2 s_1 T &=\begin{matrix}
1&4&4&8&8\\
2&5&5&9& \\
3&6&6& & \\
7&7& & & \\
8&9& & & \\
9& & & & 
\end{matrix}
\end{split}
\end{equation*}
We now give the statistics $d_{1,u R}$, $d_{2,u R}$, and
$d_{u R}$ for each of the tableaux $u T$, with $t=3$.
\begin{equation*}
\begin{matrix}
u&d_{1,u R}&d_{2,u R}&d_{u R}\\
\mathrm{id}&3&1&7\\
s_1&3&0&6\\
s_2&2&1&5\\
s_2 s_1&3&0&6\\
s_1 s_2&2&1&5\\
s_1 s_2 s_1&3&1&7
\end{matrix}
\end{equation*}
So $\charge_R(T)=(7+6+5+6+5+7)/3!=6$.  
\end{ex}

Finally we give the symmetry bijection for catabolizable tableaux.
Let $R$ and $R'$ be sequences of rectangles that
rearrange each other and $u\in\Sigma_{[n]}$ the shortest permutation
such that $\gamma(R') = u \gamma(R)$.
From \cite{BK} it follows that the following diagram commutes,
where $u$ acts by an automorphism of conjugation:
\begin{equation} \label{auto}
\begin{CD}
	CST(\la,\gamma(R)) @>{u}>> CST(\la,\gamma(R')) \\
	@V{\Phi_{\rows(R)}}VV		@VV{\Phi_{\rows(R')}}V	\\
	RC(\la^t;\cols(R^t)) @= RC(\la^t;\cols((R')^t)) 
\end{CD}
\end{equation}

Conjecture \ref{image conj} implies the following result.

\begin{conj} \label{CT symmetry}
The automorphism of conjugation $u$
restricts to a bijection 
\begin{equation*}
\Image b_R \rightarrow \Image b_{R'},
\end{equation*}
so that the diagram commutes:
\begin{equation*}
\begin{CD}
	RC(\la^t;R^t) @= RC(\la^t;(R')^t) \\
	@V{b_R}VV		@VV{b_{R'}}V	\\
	CST(\la,\gamma(R)) @>>u> CST(\la,\gamma(R'))
\end{CD}
\end{equation*}
In particular, if both $R$ and $R'$ are dominant,
then $u$ restricts to a bijection
\begin{equation*}
  CT(\la;R)\rightarrow CT(\la;R').
\end{equation*}
\end{conj}

The conjecture is trivial when each rectangle in $R$ is a single
row, but is quite interesting even when each rectangle is a single
column (in which case the permutation $u$ is the identity).

\begin{ex} Let $R'=s_1 R=((2,2,2,2),(3,3),(1,1,1))$.
Then the permutation $u$ of minimal length sending
$\gamma(R)$ to $\gamma(R')$
is given by $u=561234789$, which has the reduced word
$43215432$.  Let $S$ be the first tableau in Example \ref{cat ex}.
Here the operators $s_i$ are the automorphisms of conjugation,
or equivalently the rectangle-switching bijections for the
appropriate sequences of one-rowed rectangles.
\begin{equation*}
\begin{split}
S=\begin{matrix}
1&1&1&5&6\\
2&2&2&6& \\
3&3&7& & \\
4&4& & & \\
5&8& & & \\
9& & & &
\end{matrix}\qquad
s_4 s_3 s_2 s_1 s_5 s_4 s_3 s_2 S = \begin{matrix}
1&1&5&5&6\\
2&2&6&6& \\
3&3&7& & \\
4&4& & & \\
5&8& & & \\
9& & & &
\end{matrix}
\end{split}
\end{equation*}
\end{ex}

\section{Duality bijections}
\label{duality section}

This section uses the notation of subsection \ref{duality}.
Fix an integer $k$ such that $k\ge\la_1$ and $k\ge\mu_i$
for all $1\le i\le t$.  Let $(\Rt)^t$ denote the sequence
of rectangles whose $i$-th partition is $(\Rt_i)^t$,
a $(k-\mu_i)\times \eta_i$ rectangle.

First we define a duality bijection for rigged configurations.

\begin{prop} \label{RC dual}
There is a bijection of admissible configurations 
\begin{equation*}
\begin{split}
C(\la^t;R^t)&\rightarrow C((\lt)^t;(\Rt)^t) \\
\nu=(\nu^1,\nu^2,\dots,\nu^{k-1},(),\dots) &\mapsto
\nt=(\nu^{k-1},\nu^{k-2},\dots,\nu^1,(),\dots)
\end{split}
\end{equation*}
Furthermore, for every $1\le i \le k-1$ and $j\ge 1$,
$P_{i,j}(\nu)=P_{k-i,j}(\nt)$.
In particular the above bijection on configurations induces a
cocharge-preserving map of rigged configurations
$RC(\la^t;R^t)\rightarrow RC((\lt)^t;(\Rt)^t)$ such that
$(\nu,L)\mapsto(\nt,\Lt)$ where $\Lt^p_s=\Lt^{k-p}_s$ for all
$1\le p < k$ and $1\le s\le \ell(\nu^p)=\ell(\nt^{k-p})$.
\end{prop}

In other words, the bijection merely replaces the first $k-1$
labelled partitions of $(\nu,L)$ (the rest are empty) by the
reverse sequence of labelled partitions.  The proof is straightforward.
It follows immediately that
\begin{equation*}
  RC_{\la^t;R^t}(q) = RC_{(\lt)^t;(\Rt)^t}(q)
\end{equation*}

Next we consider the duality bijections for LR tableaux.
Let $T$ be any column strict tableau of shape $\la$ and
content $\gamma$.  Define the dual tableau $\widetilde{T}$ of $T$
(with respect to the $n\times k$ rectangle) to be the unique column strict
tableau of shape $\lt$ such that the $j^{th}$ column of $\widetilde{T}$
is the set complement within the interval $[n]$, of the $(k+1-j)$-th column
of $T$.  Clearly $\widetilde{T}$ has content
$\widetilde{\gamma}=(k-\gamma_n,k-\gamma_{n-1},\dots,k-\gamma_1)$.

\begin{ex} In our examples, we have $n=9$ and
\begin{equation*}
\la=(5,4,3,2,2,1,0,0,0) \qquad R=((3,3),(2,2,2,2),(1,1,1))
\end{equation*}
 Let $k=5$.  Then
\begin{equation*}
\begin{split}
\widetilde{\la}&=(5,5,5,4,3,3,2,1,0) \\
\widetilde{R}&=((2,2),(3,3,3,3),(4,4,4))
\end{split}
\end{equation*}
The tableau $T\in LRT(\la;R)$ is sent to $\widetilde{T}\in LRT(\lt;\Rt)$.
\begin{equation*}
T=\begin{matrix}
1&1&1&3&3\\
2&2&2&4& \\
4&5&7& & \\
5&6& & & \\
6&8& & & \\
9& & & &
\end{matrix} \qquad
\widetilde{T}=\begin{matrix}
1&1&3&3&3\\
2&2&4&4&7\\
4&5&5&7&8\\
5&6&6&9& \\
6&7&8& & \\
7&8&9& & \\
8&9& & & \\
9& & & &
\end{matrix}
\end{equation*}
\end{ex}

The duality map respects Knuth equivalence in the following sense.
Let $b$ be a strictly decreasing word in the alphabet $[n]$.
Let $\widetilde{b}$ be the strictly decreasing word whose letters
are complementary in $[n]$ to those of $b$.

\begin{prop} \cite{RS}
\label{jeu and dual} Let $a$, $b$, $c$, and $d$ be strictly
decreasing words in the alphabet $[n]$.  Then
$ab \Kn cd$ if and only if
$\widetilde{b}\widetilde{a}\Kn \widetilde{d}\widetilde{c}$
where $\Kn$ denotes Knuth equivalence.
\end{prop}

The duality bijection on column strict tableaux also
restricts to a map from LR tableaux to LR tableaux.

\begin{prop} The bijection $T\mapsto \widetilde{T}$ restricts to a
bijection 
\begin{equation*}
  LRT(\la;R)\rightarrow LRT(\lt;\Rt).
\end{equation*}
Furthermore,
for every $T\in LRT(\la;R)$ we have
\begin{equation} \label{duality and charge_R}
  \charge_R(T) = \charge_{\Rt}(\widetilde{T}),
\end{equation}
using \eqref{charge_R} as the definition of $\charge_R$.
\end{prop}
\begin{proof} Let $T$ be a column strict tableau of shape
$\la$ and content $\gamma$.  The following are equivalent:
\begin{enumerate}
\item $T\in LRT(\la;R)$.
\item $T|_{A_i}$ is lattice in the alphabet $A_i$ for all $i$.
\item For all $0\le j\le k$ and all $i$,
the last $j$ columns of $T|_{A_i}$ have partition content
in the alphabet $A_i$.
\item For all $0\le j\le k$ and all $i$,
the first $j$ columns of $\widetilde{T}_{A_i}$ has antipartition content
in the alphabet $A_i$.
\item For all $0\le j\le k$ and all $i$,
the last $k-j$ columns of $\widetilde{T}|_{A_i}$ have partition content
in the alphabet $A_i$.
\item $\widetilde{T}|_{A_i}$ is lattice in the alphabet $A_i$
for all $i$.
\item $\widetilde{T}\in LRT(\lt;\Rt)$.
\end{enumerate}
The equivalence of the first three and last three assertions follow
by definition.  The equivalence of items 3 and 4 follows from the
definition of $\widetilde{T}$.  The equivalence of 4 and 5 follows
from the fact that for each letter $x\in A_i$, the letter $x$ appears
with total multiplicity $\mu_i$.

To verify \eqref{duality and charge_R}, it suffices to show
that for all $1\le p \le t-1$ and all $u$,
\begin{equation*}
  d_{p,u R}(u T) = d_{p,u \Rt}(u \widetilde{T})
\end{equation*}
It is immediate that $u \Rt = \widetilde{u R}$.
It suffices to establish the identities
\begin{equation}
\begin{align}
\label{d and dual} d_{p,R}(T) &= d_{p,\Rt}(\widetilde{T})\qquad
		\text{for all $R$} \\
\label{perm and dual} u \widetilde{T}&=\widetilde{u T}\qquad
		\text{for all $u$}.
\end{align}
\end{equation}
To prove \eqref{d and dual}, one can reduce to the case
$p=1$ and $t=2$ by Lemma \ref{jeu and dual}.
Without loss of generality assume that $\mu_1\ge\mu_2$.  Let
$\alpha$ and $\beta$ be the left and right partitions given by slicing
the skew shape $\la/R_1$ vertically just after the $\mu_1$-th column.
The Littlewood-Richardson rule implies that the
180 degree rotation of $\alpha$ fits together with $\beta$ to
form the rectangle $R_2$.  We have
\begin{equation*}
\begin{split}
  d_{1,R}(T) &= \mu_2 n - (\la^t_1+\dots+\la^t_{\mu_2}) \\
    &= \lt^t_k +\lt^t_{k-1}+\dots+\lt^t_{k+1-\mu_2} \\
    &= d_{1,\Rt}(\widetilde{T}),
\end{split}
\end{equation*}
proving \eqref{d and dual}.

It suffices to prove \eqref{perm and dual} when $u$ is an adjacent
transposition $s_p$.  Let $B=A_p\cup A_{p+1}$.  By the definitions it
follows immediately that $s_p \widetilde{T}$ and $\widetilde{s_p T}$ agree
when restricted to the complement of the alphabet $B$.  It remains
to show that the restrictions of the two tableaux to $B$ agree.
By Lemma \ref{jeu and dual}, we may assume that $p=1$ and $t=2$.
But this case follows immediately since taking the dual tableau
and applying $s_1$ both send LR tableaux to LR tableaux, and
all of the relevant sets of LR tableaux are singletons.
\end{proof}

It follows from \eqref{charge_R} that
\begin{equation*}
  LRT_{\la;R}(q) = LRT_{\lt;\Rt}(q)
\end{equation*}

Moreover,

\begin{conj} The following diagram commutes
\begin{equation*}
\begin{CD}
	LRT(\la;R) @>\textrm{dual}>> LRT(\lt;\Rt) \\
	@V{\Phi_R}VV		@VV{\Phi_{\Rt}}V	\\
	RC(\la^t;R^t) @>>> RC(\la^t;(\Rt)^t) 
\end{CD}
\end{equation*}
where the bottom map is given in Prop. \ref{RC dual}.
\end{conj}

To discuss the duality bijection for catabolizable tableaux,
we give a modified definition that slices the tableaux vertically
rather than horizontally.  In this section let us refer to
$R$-catabolizability as $R$-row catabolizability and
to the $R_1$-catabolism as the $R_1$-row catabolism.

Let $S$ be a column strict tableau of partition shape with
$S|_{A_1}=K_1$.  Let $S_l$ and $S_r$ be the left and right
subtableaux obtained by slicing the skew tableau $S-K_1$
vertically just after the $\mu_1$-th column.
Define the $R_1$-column catabolism $\ccat_{R_1}(S)$ of $S$ by
$P(S_r S_l)$.  Say that $S$ is \textit{$R$-column catabolizable} if
if $S|_{A_1}=K_1$ and $\ccat_{R_1}(S)$ is
$\Rh$-column catabolizable, where $\Rh=(R_2,R_3,\dotsc)$.
Let $CCT(\la;R)$ denote the set of $R$-column catabolizable tableaux
of shape $\la$.

\begin{prop} \label{cat dual}
$S$ is $R$-column catabolizable if and only if
$\widetilde{S}$ is $\Rt$-column catabolizable.
\end{prop}
\begin{proof} Fix a sufficiently large number $k$.
All dual tableaux will be taken with respect to $k$ columns.
Note first that $S|_{A_1}=K_1$ if and only if
$\widetilde{S}|_{A_1}$ is the tableau of shape $\Rt_1$
whose $i$-th row consists of $k-\mu_1$ copies of the letter $i$
for all $1\le i\le \eta_1$.  Then
Proposition \ref{jeu and dual} implies that
the dual of the tableau $\ccat_{R_1}(S)$ with respect to
the alphabet $[n]-A_1$, is equal to the tableau
$\ccat_{\Rt_1}(\widetilde{S})$.  The result follows by induction.
\end{proof}

The following conjecture connects the two kinds of catabolizability
using the images of the maps $b_R$.

\begin{conj} \label{col image conj}
Suppose $R$ is a sequence of rectangles such that
$R^t$ is dominant.  Then $\Image b_R = CCT(\la;R)$.
\end{conj}

Using Conjectures \ref{CT symmetry} and \ref{col image conj}
one obtains a duality bijection for row catabolizable tableaux
using automorphisms of conjugation and the tableau duality
bijection.

\begin{ex}
If the hypothesis of Conjecture \ref{col image conj} is not satisfied
then $CCT(\la;R)$ could be too large. For example,
let $\la=(2,2)$ and $R=((1,1,1),(1))$. Then $CCT(\la;R)$ is empty.
Let $R'=((1),(1,1,1))$. $CCT(\la;R')$ is not empty; in fact it is
equal to $CCT((2,2);((1,1),(1,1)))$.
\end{ex}

\section{Transpose Bijections}
\label{transpose section}

A bijection $LRT(\lambda;R)\rightarrow LRT(\lambda^t;R^t)$
is given by the relabelling which sends $T$ to the 
transpose of the tableau obtained from $T$
by replacing the $j$-th occurrence (from the left) of the letter
$\eta_1+\eta_2+\dots+\eta_{i-1}+k$ by the letter $\mu_1+\mu_2+\dots+
\mu_{i-1}+j$, for all $i$, $1\le j\le \mu_i$ and $1\le k \le \eta_i$.
We call this map the LR transpose.

\begin{ex} $\la^t=(6,5,3,2,1)$ and $R_1^t=(2,2,2)$,
$R_2^t=(4,4)$ and $R_3^t=(3)$.
The set $LRT(\la^t;R^t)$ is given by the following
four tableaux, which are the images under the LR transpose
map of the four tableaux of $LRT(\la;R)$ listed in
Example \ref{LR example}.
\begin{equation*}
\begin{matrix}
1&1&4&4&4&6\\
2&2&5&5&6& \\
3&3&6& & & \\
4&5& & & & \\
5& & & & &
\end{matrix}\quad
\begin{matrix}
1&1&4&4&4&6\\
2&2&5&5&5& \\
3&3&6& & & \\
5&5& & & & \\
6& & & & &
\end{matrix}\quad
\begin{matrix}
1&1&4&4&4&4\\
2&2&5&5&6& \\
3&3&6& & & \\
5&5& & & & \\
6& & & & &
\end{matrix}\quad
\begin{matrix}
1&1&4&4&4&4\\
2&2&5&5&5& \\
3&3&6& & & \\
5&6& & & & \\
6& & & & &
\end{matrix}
\end{equation*}
\end{ex}

For rigged configurations, 
we wish to define a bijection $RC(\la;R)\rightarrow RC(\la^t;R^t)$
(called the RC-transpose)
sending $(\nu,L)\rightarrow (\nh,\Lh)$ with the property that
\begin{equation} \label{stat comp}
  \cocharge(\nu,L) = n(R) - \cocharge(\nh,\Lh).
\end{equation}

The bijection on configurations is most easily defined
using a variant of the original construction.
Let $\nu$ be an admissible $(\la;R)$ configuration.  Recall that
$\alpha_{k,n}$ is the $n$-th part of the transpose of the partition $\nu^k$
and that $\nu^0$ is the empty partition.
Define the matrix $(m_{i,j})$ by
\begin{equation} \label{m def}
	m_{i,j} = \alpha_{i-1,j}-\alpha_{i,j}
\end{equation}
for $i,j\ge 1$, or equivalently,
\begin{equation} \label{inv m def}
	\alpha_{i,j} = - \sum_{k=1}^i m_{k,j}
\end{equation}
The matrix $(m_{i,j})$ will be used in place of the configuration $\nu$.
Before defining the transpose map on rigged configurations, let us
calculate the row and column sums of the matrix $(m_{i,j})$ and
the cocharge of $\nu$ in terms of the $m_{i,j}$.
Let $\theta$ be the indicator function for nonnegative numbers:
\begin{equation*}
\theta(x)=
  \begin{cases} 1 & \text{ if $x\ge 0$,} \\
		0 & \text{ if $x < 0$.}
  \end{cases}
\end{equation*}
Since $\alpha_{0,j}=0$ for all $j$ and $\alpha_{i,j}=0$ for large $i$,
\begin{equation}
\label{m col sum}
	\sum_i m_{i,j} = \sum_i (\alpha_{i-1,j}-\alpha_{i,j}) = 0.
\end{equation}
Using the definition of a configuration of type $(\la;R)$
and the notation $r_+ = \max(r,0)$, we have
\begin{equation}
\label{m row sum}
\begin{split}
	\sum_j m_{i,j} &= \sum_j (\alpha_{i-1,j}-\alpha_{i,j}) \\
	&= |\nu^{i-1}| - |\nu^i| \\
	&= \sum_{j>i-1} (\la_j - \sum_a \mu_a (\eta_a-(i-1))_+) -
	\sum_{j>i} (\la_j - \sum_a \mu_a (\eta_a - i)_+) \\
	&= \la_i - \sum_a \mu_a \theta(\eta_a-i),
\end{split}
\end{equation}
The cocharge \eqref{RC cocharge} can be rewritten as
\begin{equation*}
\begin{split}
  \cocharge(\nu)
  &= \sum_{k,n\ge 1} \alpha_{k,n}(\alpha_{k,n}-\alpha_{k+1,n}) \\
  &= \sum_{k,n\ge 1} (-m_{1,n}-m_{2,n}-\dots-m_{k,n}) m_{k+1,n} \\
  &= - \sum_{\substack{n\ge 1 \\ j>k\ge 1}} m_{k,n} m_{j,n} \\
  &= - 1/2 \sum_{\substack{n \ge 1 \\ j \not= k\ge 1}} m_{k,n} m_{j,n} \\
  &= 1/2 \sum_{k,n\ge 1} m_{k,n}^2 - 1/2 \sum_{j,k,n\ge 1} m_{k,n} m_{j,n} \\
  &= 1/2 \sum_{k,n\ge 1} (m_{k,n}^2-m_{k,n}) + 0 \\
  &= \sum_{k,n\ge 1} \binom{m_{k,n}}{2} =\\
  &= 1/2 \sum_{k,n\ge 1} m_{k,n}^2.
\end{split}
\end{equation*}

Define the matrix $\mh_{i,j}$ by
\begin{equation}
\label{m trans def}
	\mh_{i,j} = - m_{j,i} + \theta(\la_j - i) -
	\sum_a \theta(\mu_a - i)\theta(\eta_a - j)
\end{equation}
Observe that this process is an involution.  Applying
this process to $(\mh_{i,j})$ with respect to the pair $(\la^t;R^t)$,
\begin{equation*}
\begin{split}
\widehat{\mh}_{i,j} &= - \mh_{j,i} + \theta(\la^t_j-i)-
	\sum_a \theta(\eta_a-i)\theta(\mu_a-j) %\\
\end{split}
\end{equation*}
\begin{equation*}
\begin{split}
  &= -(- m_{i,j} + \theta(\la_i-j)-\sum_a \theta(\mu_a-j)\theta(\eta_-i))\\
	&+\theta(\la^t_j-i)-\sum_a \theta(\eta_a-i)\theta(\mu_a-j) \\
  &= m_{i,j} +\theta(\la^t_j-i)-\theta(\la_i-j) = m_{i,j}
\end{split}
\end{equation*} 

Note that in \eqref{m trans def} the sum over $a$ is equal to $r_{j,i}(R)$,
the number of rectangles in $R$ that contain the cell $(j,i)$.
We must show that the matrix $(\mh_{i,j})$ corresponds to an
admissible configuration of type $(\la^t;R^t)$.
Let $\alh_{i,j}$, $\nh^i_j$, $P_{i,j}(\nh)$, and $m_i(\nh^j)$ denote
the analogous quantities involving $\mh_{i,j}$ in place of $m_{i,j}$.  
The first step is to show that each $\nh^i$ is a partition,
that is, $m_j(\nh^i)\ge 0$ for all $i,j\ge 1$.  We have
\begin{equation} \label{mt pos}
\begin{split}
  m_j(\nh^i) &= \alh_{i,j} - \alh_{i,j+1} \\
  &= \sum_{k=1}^i (-\mh_{k,j}+\mh_{k,j+1}) \\
  &= \sum_{k=1}^i (m_{j,k}-m_{j+1,k}-\theta(\la_j-k)+\theta(\la_{j+1}-k)\\
	&+\sum_a \theta(\mu_a-k)(\theta(\eta_a-j)-\theta(\eta_a-(j+1))) \\
  &= \sum_{k=1}^i (\alpha_{j-1,k}-\alpha_{j,k}-\alpha_{j,k}+\alpha_{j+1,k})\\
  	& -\min(\la_j,i)+\min(\la_{j+1},i)
  	+ \sum_a \min(\mu_a,i) \delta_{\eta_a,j} \\
  &= Q_i(\nu^{j-1})-2 Q_i(\nu^j)+Q_i(\nu^{j+1}) \\
  	&+ \sum_a \min(\mu_a,i) \delta_{\eta_a,j}
	-\min(\la_j,i)+\min(\la_{j+1},i) \\
  &= P_{j,i}(\nu) - \min(\la_j,i)+\min(\la_{j+1},i)
\end{split}
\end{equation}

We require the following technical result on vacancy numbers, whose
proof is in the appendix.  Recall that
$m_n(\rho)$ is the number of parts of the partition $\rho$ of size $n$.

\begin{lem} \label{vacancy} Let $\nu$ be a configuration of type $(\la;R)$.
The following are equivalent.
\begin{enumerate}
\item $\nu$ is admissible, that is, $P_{k,n}(\nu)\ge 0$ for all $k,n\ge 1$.
\item For every $k,n\ge 1$, if $m_n(\nu^k)>0$ then $P_{k,n}(\nu)\ge 0$.
\item For every $k,n\ge 1$,
\begin{equation} \label{more vacancy}
  P_{k,n}(\nu)\ge \min(\la_k,n)-\min(\la_{k+1},n)
\end{equation}
\end{enumerate}
Moreover, if $\nu$ is admissible then
\begin{equation}
  m_n(\nu^k) = 0 \qquad\text{ whenever $n>\la_{k+1}$.}
\end{equation}
\end{lem}

From \eqref{mt pos}, Lemma \ref{vacancy}, and the admissibility of
$\nu$, it follows that $\nh^i$ is a partition for all $i\ge 1$.
Next it is verified that $\nh$ is a configuration of type $(\la^t;R^t)$.
\begin{equation*}
\begin{split}
  |\nh^i| &= \sum_j \alh_{i,j} = \sum_j \sum_{k=1}^i -\mh_{k,j} \\
  &= \sum_j \sum_{k=1}^i (m_{j,k}-\theta(\la_j-k)
  	+\sum_a \theta(\mu_a-k)\theta(\eta_a-j)) \\
  &= \sum_{k=1}^i \sum_j m_{j,k} - \sum_j \min(\la_j,i) +
  \sum_a \min(\mu_a,i) \eta_a \\
  &= - \sum_{j=1}^i \la^t_j + \sum_a \min(\mu_a,i) \eta_a
\end{split}
\end{equation*}
by \eqref{m trans def} and \eqref{m col sum}.  Comparing this with
\eqref{config}, $\nh$ is a configuration of type $(\la^t;R^t)$,
since $R^t$ is obtained from $R$ by switching the roles of $\mu$ and $\eta$.
Finally it must be verified that $\nh$ is admissible.
Since the map $m\mapsto\mh$ is involutive as a map of matrices, it is valid
to apply the formula \eqref{mt pos} to $\nh$ instead of $\nu$, obtaining
\begin{equation*}
  P_{i,j}(\nh) = m_i(\nu^j) + \min(\la^t_i,j)-\min(\la^t_{i+1},j)
	\ge m_i(\nu^j) \ge 0.
\end{equation*}
Therefore the map $(m_{i,j})\mapsto(\mh_{i,j})$ defines a bijection
$C(\la;R)\rightarrow C(\la^t;R^t)$.

This map is extended to riggings as follows.
By Lemma \ref{vacancy} and \eqref{mt pos}, the map $m\mapsto \mh$ 
has the additional property that 
\begin{equation} \label{trans m and P}
\begin{split}
  m_n(\nu^k) &= P_{n,k}(\nh) \\
  P_{k,n}(\nu) &= m_k(\nh^n)
\end{split}
\end{equation}
provided that if $m_n(\nu^k)>0$.
Applying this to the inverse map, if $m_k(\nt^n)>0$ then
these same two equalities hold.  This implies that the two sets of
rigged configurations $RC(\la;R)$ and $RC(\la^t;R^t)$ have the same
cardinality.  Let us give a specific bijection between the riggings.

Let $(\nu,L)\in RC(\la;R)$.  Let $\nh$ be the admissible
configuration of type $(\la^t;R^t)$ given by \eqref{m trans def}.
Note that a rigging $L$ of $\nu$ determines, for each pair
$k,n\ge 1$, a partition $\rho_{k,n}(\nu,L)$ inside a rectangle
of height $m_n(\nu^k)$ and width $P_{k,n}(\nu)$ given by the
the labels of the parts of $\nu^k$ of size $n$.

Suppose $m_n(\nu^k)>0$.  Let $\Lh$ be the rigging
of the configuration $\nh$ such that
$\rho_{n,k}(\nh,\Lh)$ is the transpose of the
complementary partition to $\rho_{k,n}(\nu,L)$ in the
rectangle of height $m_n(\nu^k)$ and width $P_{k,n}(\nu)$ rectangle,
for all $k,n\ge1$.

Then the map $(\nu,L)\mapsto(\nh,\Lh)$ defines the RC-transpose bijection
$RC(\la;R)\rightarrow RC(\la^t;R^t)$.

\begin{prop} \label{transpose poly}
\begin{equation}
	RC_{\la^t;R^t}(q) = q^{n(R)} RC_{\la;R}(q^{-1})
\end{equation}
\end{prop}
\begin{proof} It is enough to check that the RC-transpose bijection
$(\nu,L)\mapsto(\nh,\Lh)$ satisfies
\begin{equation} \label{transpose cocharge}
\cocharge(\nu,L)+\cocharge(\nh,\Lh)=n(R)
\end{equation}
(see \eqref{n statistic}).
By the definition of the rigging $\Lh$, it is enough to check that
\begin{equation} \label{costatistic}
  \cocharge(\nu) + \cocharge(\nh) + \sum_{k,n\ge1} P_{k,n}(\nu) m_n(\nu^k)
  = n(R)
\end{equation}

The sum $\sum_{k,n} P_{k,n}(\nu)m_n(\nu^k)$ is calculated first.
For the vacancy numbers,
\begin{equation*}
\begin{split}
  Q_n(\nu^{k-1})-2 Q_n(\nu^k)+Q_n(\nu^{k+1})
  &= \sum_{j=1}^n (\alpha_{k-1,j}-2 \alpha_{k,j}+\alpha_{k+1,j}) \\
  &= \sum_{j=1}^n (m_{k,j}-m_{k+1,j})
\end{split}
\end{equation*}
and
\begin{equation} \label{P in m}
P_{k,n}(\nu) = \sum_{j=1}^n (m_{k,j}-m_{k+1,j})
	+ \sum_{a\ge1} \min(\mu_a,n) \delta_{\eta_a,k}.
\end{equation}
The multiplicities $m_n(\nu^k)$ are given by
\begin{equation} \label{m in m}
  m_n(\nu^k) = \alpha_{k,n}-\alpha_{k,n+1} =
	\sum_{i=1}^k (-m_{i,n}+m_{i,n+1}).
\end{equation}
The desired sum is given by
\begin{equation} \label{P times m}
\begin{split}
  \sum_{k,n\ge1} P_{k,n}(\nu)m_n(\nu^k) &=
  \sum_{k,n\ge1} (\sum_{j=1}^n (m_{k,j}-m_{k+1,j}) \\
	&+	  \sum_a \min(\mu_a,n)\delta_{\eta_a,k})
		 \sum_{i=1}^k (-m_{i,n}+m_{i,n+1}) \\
  &= \sum_{\substack{k\ge i\ge1 \\ n\ge j\ge1}}
    (m_{k,j}-m_{k+1,j})(-m_{i,n}+m_{i,n+1})  \\
  & + \sum_{\substack{k\ge i\ge 1 \\ n\ge 1}} 
  \sum_a \min(\mu_a,n)\delta_{\eta_a,k}) (-m_{i,n}+m_{i,n+1}) \\
  &= - \sum_{i,j\ge 1} m_{i,j}^2 +
  \sum_a \sum_{i=1}^{\eta_a}
	\sum_{n\ge 1} \min(\mu_a,n)(-m_{i,n}-m_{i,n+1}) \\
  &= - 2\, \cocharge(\nu)
     - \sum_a \sum_{i=1}^{\eta_a} \sum_{n=1}^{\mu_a} m_{i,n} \\
  &= - 2\, \cocharge(\nu)
     - \sum_a \sum_{i,j\ge1} m_{i,j}\theta(\mu_a-j)\theta(\eta_a-i) \\
  &= - 2\, \cocharge(\nu) - \sum_{i,j\ge1} m_{i,j} r_{i,j}
\end{split}
\end{equation}
where $r_{i,j}=r_{i,j}(R)$ is as in the definition of $n(R)$.
\eqref{n statistic}

Next let us calculate the cocharge of the configuration $\nh$.
By \eqref{RC cocharge} and \eqref{m trans def}, we have
\begin{equation}
  \cocharge(\nh) = \sum_{i,j\ge1} \binom{\mh_{i,j}}{2} 
  = \binom{-m_{j,i}+\theta(\la_j-i) - r_{j,i}}{2}.
\end{equation}
Using the identities $\binom{a+b}{2}=\binom{a}{2}+\binom{b}{2}+ab$
and $\binom{-a}{2}=\binom{a}{2}+a$, we have
\begin{equation*}
\begin{split}
  \cocharge(\nh) &= \sum_{i,j\ge1} (
    \binom{m_{j,i}}{2}+\binom{\theta(\la_j-i)}{2}+\binom{r_{j,i}}{2} \\
  &  -m_{j,i}\theta(\la_j-i)+m_{j,i}r_{j,i}-\theta(\la_j-i)r_{j,i}
	+m_{j,i} + r_{j,i})\\
  &= \cocharge(\nu) + 0 + n(R) +
	\sum_{i,j\ge1} (m_{j,i}r_{j,i} + (1-\theta(\la_j-i))
	(m_{j,i}+r_{j,i}).
\end{split}
\end{equation*}
Now $\theta(\la_j-i)=0$ if and only if the cell $(j,i)$ is not
in the partition shape $\la$.  Since the number of rigged configurations
$|RC(\la;R)|$ is equal to the LR coefficient $\LRC^\la_R$, it follows that
each rectangle $R_a$ is contained in $\la$.  Thus $\theta(\la_j-i)=0$
implies that $r_{j,i}=0$.  Lemma \ref{vacancy} guarantees that
$m_{j,i}=\alpha_{j-1,i}-\alpha_{j,i}=0$ if $\theta(\la_j-i)=$.  Thus
\begin{equation*}
  \cocharge(\nh) = \cocharge(\nu) + n(R) + \sum_{i,j\ge1} m_{j,i}r_{j,i}
\end{equation*}
Together with \eqref{P times m}, this implies \eqref{costatistic}.
\end{proof}

The transpose bijections were chosen with the following property
in mind.

\begin{conj} \label{LR RC trans comm}
The following diagram commutes:
\begin{equation*}
\begin{CD}
	LRT(\la;R) @>{\mathrm{LR-transpose}}>> LRT(\la^t;R^t) \\
	@V{\Phi_R}VV		@VV{\Phi_{R^t}}V	\\
	RC(\la^t;R^t) @>{\mathrm{RC-transpose}}>> RC(\la;R) 
\end{CD}
\end{equation*}
\end{conj}

We now give a map from the $R$-catabolizable tableaux of shape $\la$,
to the $R$-column catabolizable tableaux of shape $\lt$, in the
case that $R$ is dominant.

Let $S\in CT(\la;R)$.  Let $S_1=\cat_{R_1}(S)$ and let
$Q_1=Q_c(\word(S_+)\word(S_-))$ where $\word(T)$ is the
row-reading word of the (skew) tableau $T$, and
$Q_c(w)$ is the Schensted $Q$ symbol for the column insertion of the word
$w$ \cite{Sch}.  The image $U\in CT(\lambda^t;R^t)$ of $S$
under our proposed bijection, is uniquely defined by the property
that $\ccat(U)=U_1$ where $U_1$ is the image of $S_1$ (which
is defined by induction), and that $Q_1^t=Q(\word(U_r)\word(U_l))$,
where $Q(a)$ is the Schensted $Q$ symbol for the row insertion
of $a$  (see the definition of $\ccat_{R_1}(U)$).

In practice one first computes the sequence of tableaux $S_0=S$ and
$S_i=\cat_{R_i}(S_{i-1})$, together with the recording tableaux $Q_i$
coming from the column insertion of the appropriate row-reading words.
Then one computes the sequence of tableaux $\dotsc,U_2,U_1$ by
setting $U_i$ to be the tableau such that $\ccat_{R_i^t}(U_i)=U_{i+1}$
with recording tableau $Q_i^t$.

\begin{ex} Let $S$ be the first tableau listed in Example \ref{cat ex}.
We give the successive catabolisms of $S$
together with the recording tableaux $Q_i$.
\begin{equation}
\begin{split}
S&=S_0=
\begin{matrix}
1&1&1&5&6\\
2&2&2&6& \\
3&3&7& & \\
4&4& & & \\
5&8& & & \\
9& & & &
\end{matrix} \\
S_1&=
\begin{matrix}
3&3&7\\
4&4&8\\
5&5&9\\
6&6& 
\end{matrix}
\qquad
Q_1=
\begin{matrix}
1&2&3\\
4&5&10\\
6&7&11\\
8&9& \\
11& & 
\end{matrix} \\
S_2&=
\begin{matrix}
7\\8\\9
\end{matrix}
\qquad
Q_2=
\begin{matrix}
1\\2\\3
\end{matrix} \\
S_3 &= \emptytab
\qquad
Q_3 = \emptytab
\end{split}
\end{equation}

The sequence of tableaux $U_i$ are calculated by
``reverse column catabolisms" whose row insertions are
recorded by the transposes of the $Q$s.

\begin{equation}
\begin{split}
U_3 &= \emptytab \qquad
Q_3^t = \emptytab \\
U_2 &=
\begin{matrix}
6&6&6
\end{matrix}
\qquad
Q_2^t =
\begin{matrix}
1&2&3
\end{matrix} \\
U_1 &=
\begin{matrix}
4&4&4&4\\
5&5&5&5\\
6&6&6&
\end{matrix}
\qquad
Q_1^t =
\begin{matrix}
1&4&6&8\\
2&5&7&9\\
3&10&11& 
\end{matrix} \\
U&=U_0=
\begin{matrix}
1&1&4&4&5&6\\
2&2&5&5&6& \\
3&3&6& & & \\
4&4& & & & \\
5& & & & & 
\end{matrix}
\end{split}
\end{equation}
\end{ex}

\begin{conj} \label{CT transpose}
The above map gives a bijection $CT(\la;R)$ to
$CCT(\la^t;R^t)$ when $R$ and $R^t$ are dominant.
\end{conj}

The essential point to check is that the tableau $U$ constructed
by this algorithm is a column strict tableau, since it is conceivable
that there could be violations of column strictness between pairs
of entries of the form $U(i,\mu_1)$ and $U(i,\mu_1+1)$ for $i>\eta_1$.

\section{Monotonicity maps}
\label{monotonicity section}

For this section assume that  $R\dom R'$.

From the definitions it follows directly
that $RC(\la^t;R^t)\subseteq RC(\la^t;(R')^t)$
which is obviously a cocharge-preserving embedding.
At the end of this section we give a direct description of a related
charge-preserving embedding $\zeta:RC(\la;R)\hookrightarrow RC(\la;\rows(R))$.

A consequence of Conjecture \ref{image conj} is that if both $R$
and $R'$ are dominant, then $CT(\la;R)\subseteq CT(\la;R')$.

For LR tableaux, we describe embeddings
\begin{equation*}
  \theta^{R'}_R: LRT(\la;R)\rightarrow LRT(\la;R')
\end{equation*}
where $R\dom R'$.  To define such embeddings it is enough to assume
that $R$ covers $R'$, that is, for some $k$,
$\tau^j(R)=\tau^j(R')$ for all $k\not=k$ and
$\tau^k(R)$ covers $\tau^k(R')$ in the dominance order.
Using the rectangle switching bijections
we may assume that $R$ and $R'$ have the form
\begin{equation*}
\begin{split}
  R &= ((k^a),(k^b),R_3,R_4,\dots) \\
  R'&= ((k^{a-1},(k^{b+1}),R_3,R_4,\dots)
\end{split}
\end{equation*}
where $a>b+1$.  The injection is defined as follows; it generalizes
the rectangle-switching bijection in Section \ref{symmetry section}.
Let $T\in LRT(\la;R)$ and $B = [a+b]$.  Then $T|_B$ comprises the set
$LRT(\rho;((k^a),(k^b)))$, where $\rho$ is the shape of $T|_B$.
It follows that $LRT(\rho;((k^{a-1},k^{b+1})))$ consists of a single
tableau $T"$.  Let $T'$ be defined by
$T'|_B = T"$ and $T'|_{[n]-B} = T|_{[n]-B}$.  It is clear from the
definitions that $T'\in LRT(\la;R')$.  The injection is given by
$T\mapsto T'$.

By composing embeddings of the form $\theta^{R'}_R$ for $R\dom R'$
a covering relation, one may obtain maps of the form $\theta^{R'}_R$
for any pair $R\dom R'$.

This conjecture is a consequence of the following.

\begin{conj} \label{LR RC embed comm} Let $R\dom R'$ be a covering
relation.  Then $\theta^{R'}_R$ is independent of the sequence
of covering relations in $\dom$ leading from $R$ to $R'$, and
\begin{equation*}
\begin{CD}
	LRT(\la;R) @>{\theta^{R'}_R}>> LRT(\la;R') \\
	@V{\Phi_R}VV		@VV{\Phi_{R'}}V	\\
	RC(\la^t;R^t) @>{\mathrm{inclusion}}>> RC(\la^t;(R')^t)
\end{CD}
\end{equation*}
\end{conj}

We conclude this section with a charge-preserving embedding 
\begin{equation*}
  \zeta~:~RC(\la;R)\hookrightarrow RC(\la;\rows(R)).
\end{equation*}

First we must define the charge of a rigged configuration
and say a few words about quantum versus coquantum numbers.

The charge of a rigged configuration $(\nu,J)$ is defined by
\begin{equation*}
\charge(\nu,J)=\charge(\nu)+\sum_{k\ge 1}\sum_{s=1}^{l(\nu^{k})}J_s^k,
\end{equation*}
where
\begin{equation*}
\charge(\nu)=\sum_{k,n\ge 1}\binom{\alpha_{k-1,n}-\alpha_{k,n}+
\sum_a \theta(\eta_a-k)\theta(\mu_a-n)}{2}.
\end{equation*}
The definition of charge is compatible with the quantum number labeling
in the following sense.  Suppose $(\nu,L)\in RC(\la;R)$; we view $L$ as
coquantum numbers.  Let $J$ be the rigging of $\nu$ obtained by
complementing the rigging $L$, that is,
\begin{equation*}
  J^k_s = P_{k,\nu^k_s}(\nu) - L^k_s.
\end{equation*}
Write $\Omega_R(\nu,L)=(\nu,J)$ for this involution.
The $J$ are quantum numbers.  The point is that
\begin{equation} \label{charge and cocharge}
  \charge(\nu,J) + \cocharge(\nu,L) = n(R).
\end{equation}
To see this, it is equivalent to show that
\begin{equation*}
  \charge(\nu) = n(R) - \cocharge(\nu) -\sum_{k,n\ge1} P_{k,n}(\nu)m_n(\nu^k).
\end{equation*}
In light of \eqref{P times m} this is equivalent to
\begin{equation*}
  \charge(\nu) = n(R) + \cocharge(\nu) + \sum_{k,n\ge1} m_{k,n} r_{k,n}.
\end{equation*}
where $(m_{k,n})$ is the matrix associated to $\nu$ (see \eqref{m def})
and $r_{k,n}=r_{k,n}(R)$ is as in \eqref{n statistic}.  Then
\begin{equation*}
\begin{split}
\charge(\nu)&=\sum_{k,n\ge 1}\binom{m_{k,n}+r_{k,n}}{2} \\
  &= \sum_{k,n\ge1} (\binom{m_{k,n}}{2}+\binom{r_{k,n}}{2}+m_{k,n}r_{k,n})\\
  &= \cocharge(\nu)+n(R)+\sum_{k,n} m_{k,n} r_{k,n}
 \end{split}
\end{equation*}
which proves \eqref{charge and cocharge}.

Now let us define the map $\zeta$.
Let $\nu=(\nu^1,\nu^2,\ldots )$ be an
admissible configuration of type $(\la;R)$ and $J$ a rigging of $\nu$.
Define $\zeta(\nu,J)=(\nt,\Jt)$ where
\begin{equation*}
  \nt^k=\nu^k\bigcup_{a\ge 1}\left(\mu_a^{(\eta_a-k)_+}\right)
\end{equation*}
and the quantum numbers on the ``old" rows remain the same,
and the ``new" rows are assigned the quantum number 0.

\begin{prop} $\zeta:RC(\la;R)\hookrightarrow RC(\la;\rows(R))$
is a charge-preserving embedding.
\end{prop}

This is an immediate consequence of the following description of $\zeta$.

\begin{lem} The map $\zeta$ satisfies the commutative diagram
\begin{equation*}
\begin{CD}
  RC(\la;R) @>{\zeta}>> RC(\la;\rows(R)) \\
    @V{\Omega_R}VV		@AA{\Omega_{\rows(R)}}A \\
  RC(\la;R) 	@.	RC(\la;\rows(R)) \\
    @V{\mathrm{transpose}_R}VV	@AA{\mathrm{transpose}_{\rows(R)}}A \\  
  RC(\la^t;R^t) @>>{{\mathrm{inclusion}}}> RC(\la^t;(\rows(R))^t)
\end{CD}
\end{equation*}
\end{lem}
\begin{proof} Let $(\nu,J)\in RC(\la;R)$ and $(\nt,\Jt)$ the image
of the composite map defined by the commutative diagram.
Let $\Omega_R(\nu,J)=(\nu,L)$, $\Omega_{\rows(R)}(\nt,\Jt)=(\nt,\Lt)$, and
let the RC-transpose of $(\nu,L)$ be $(\nh,\Lh)$.  Observe that
the composite map is charge-preserving.  We have
\begin{equation*}
\begin{split}
  \charge(\nt,\Jt) &= n(\rows(R))-\cocharge(\nt,\Lt) \\
                   &= \cocharge(\nh,\Lh) \\
                   &= n(R)-\cocharge(\nu,L) \\
		   &=\charge(\nu,J)
\end{split}
\end{equation*}
by applications of \eqref{charge and cocharge} and
\eqref{transpose cocharge} for both $R$ and $\rows(R)$.

So it suffices to show that the output $(\nt,\Jt)$ of the 
composite map agrees with $\zeta(\nu,J)$.

Let $(m_{i,j})$, $(\mh_{i,j})$, and $(\mt_{i,j})$ be the matrices
\eqref{m def} corresponding to $\nu$, $\nh$, and $\nt$.  Let
$r_{i,j}=r_{i,j}(R)$ and $\rt_{i,j}=r_{i,j}(\rows(R))$ in the notation
of \eqref{n statistic}.  By two
applications of the definition of the RC-transpose map \eqref{m trans def}
and the definition of $\rows(R)$, we have
\begin{equation*}
\begin{split}
	\mh_{i,j} &= - m_{j,i} + \theta(\la_j - i) - r_{j,i} \\
	&= - \mt_{j,i} + \theta(\la_j - i) - \rt_{j,i}.
\end{split}
\end{equation*}
Solving for $\mt_{j,i}$, we have
\begin{equation*}
  \mt_{j,i} = m_{j,i} + r_{j,i} - \rt_{j,i}.
\end{equation*}
By \eqref{inv m def} we have
\begin{equation*}
\begin{split}
  \alt_{j,i} &= \alpha_{j,i} - \sum_{k=1}^j r_{j,i}-\rt_{j,i} \\
  &= \alpha_{j,i} - \sum_{k=1}^j
	\sum_a \theta(\mu_a-i)(\theta(\eta_a-k)-\eta_a \delta_{1,k}) \\
  &= \alpha_{j,i} - \sum_a \theta(\mu_a-i) (\min(\eta_a,j)-\eta_a) \\
  &= \alpha_{j,i} + \sum_a \theta(\mu_a-i) (\eta_a-j)_+.
\end{split}
\end{equation*}
This is equivalent to the assertion that the configuration
in $\zeta(\nu,J)$ is equal to $\nt$.

To check the riggings, consider an index $k$ and a part size $n$.
It is first shown that the following two assertions suffice.
\begin{enumerate}
\item If $m_n(\nu^k)>0$ and $P_{k,n}(\nu)>0$ then
$\rho_{k,n}(\nu,J)=\rho_{k,n}(\nt,\Jt)$.
\item If $m_n(\nt^k)>0$ and $P_{k,n}(\nt)>0$ then
$\rho_{k,n}(\nu,J)=\rho_{k,n}(\nt,\Jt)$.
\end{enumerate}
To see why these conditions suffice, let us consider
a part $n$ of $\nt^k$, so that $m_n(\nt^k)>0$.
Suppose first that this part $n$ is ``new''.
If $P_{k,n}(\nt)=0$ then $\Jt$ obeys the definition of $\zeta$.
So suppose $P_{k,n}(\nt)>0$.  Then by 2,
$\rho_{k,n}(\nu,J)=\rho_{k,n}(\nt,\Jt)$.  Here the old parts of length
$n$ in $\nt^k$ have the same labels as in $\nu^k$, using up
all the nonzero parts of $\rho_{k,n}(\nt,\Jt)$.  Therefore the
new parts $n$ in $\nt^k$ have label zero, and $(\nt,\Jt)$ is
as in the definition of $\zeta$.  Otherwise suppose $n$ is an old
part of $\nt^k$, so that $m_n(\nt^k)>0$ and $m_n(\nu^k)>0$.
If either $P_{k,n}(\nu)>0$ or $P_{k,n}(\nt)>0$ then
by 1 or 2 we are done as before.  If both are zero, then
by admissibility the quantum number for this part must be
zero in both $(\nu,J)$ and $(\nt,\Jt)$, and again we are done.

Now the proof of 2 is given; the proof of 1 is similar.
For 2, the hypotheses allow us to apply \eqref{trans m and P}
for $(\nt,\Lt)\in RC(\la;\rows(R))$, so that 
$\rho_{n,k}(\nh,\Lh)$ is the complement of the transpose of
$\rho_{k,n}(\nt,\Lt)$ in the $m_n(\nt^k)\times P_{k,n}(\nt)$ rectangle.
But by the definition of $\Omega_{\rows(R)}$, 
$\rho_{k,n}(\nt,\Jt)$ is the complement of
$\rho_{k,n}(\nt,\Lt)$ in the $m_n(\nt^k)\times P_{k,n}(\nt)$ rectangle.
It follows that $\rho_{n,k}(\nh,\Lh)$ is the transpose of
the partition $\rho_{k,n}(\nt,\Jt)$, and that
$m_k(\nh^n)>0$ and $P_{n,k}(\nh)>0$.  Applying the same argument
for $(\nh,\Lh)\in RC(\la^t;R^t)$ it follows that
$\rho_{k,n}(\nu,L)$ is the transpose of the partition
$\rho_{n,k}(\nh,\Lh)$, so that 2 follows.
\end{proof}

It follows from the explicit description of $\zeta$ that
we can embed the set $RC(\la;R)$ into $RC(\la;\rows(R))$ such that for
each $k\ge 1$ and $a\ge 1$ there exists a nonnegative integer $m:=m(k,a)$
such that

$\bullet$ $\alpha_{k,m+i}-\alpha_{k,m+i+1}=\mu_a$,

$\bullet$ $J_{m+i}^{k}=0$, for all $0\le i\le (\eta_a-k)_+$.

If this characterization of the set $RC(\la;R)$ is compared with
that given by S.~Fishel \cite{Fi} for the set $M^k_{r-k}$,
it follows immediately that 
\begin{equation*}
  M^k_{r-k}(t) = \sum_{(\nu,J)\in RC(\la;R)} q^{\charge(\nu,J)}
\end{equation*}
where $R=\left((2)^{r-k},(1,1)^k,(1)^{n-2r}\right)$.
This, together with \eqref{diff poly}, implies that the
Poincar\'e polynomial $K_{\la;R}(q)$ and the above generating
function over rigged configurations, coincide in the
case that $R$ consists of some rectangles $(2)$ followed by
a sequence of rectangles $(1,1)$ and $(1)$.

\section{Appendix}

Before proving Lemma \ref{vacancy} it should be mentioned that
the implication that 2 implies 1 is fundamental.  It is used implicitly
in the definition of the bijection from LR tableaux to rigged
configurations, even in the Kostka case.  So all the results depending
upon that bijection require this lemma.  Similarly, the definition
of rigged configuration used in \cite{Fi} requires this lemma.

\begin{proof} (of Lemma \ref{vacancy}) 
Clearly it is enough to show that 2 implies 1 and 1 implies 3.
Suppose 2 holds.  It is consistent with the definitions to make the
convention that $P_{k,0}(\nu) = 0$ for all $k\ge 1$.  Using the notation
$\alpha_{k,n}$ for the $n$-th part of the transpose of the partition
$\nu^k$, the difference of vacancy numbers is given by
\begin{equation*}
  P_{k,n}(\nu) - P_{k,n-1}(\nu) =
	\alpha_{k-1,n}-2\alpha_{k,n}+\alpha_{k+1,n} +
		\sum_a \theta(\mu_a-n) \delta_{\eta_a,k},
\end{equation*}
which holds for $k,n\ge 1$.  Taking differences again,
\begin{equation*}
\begin{split}
&-P_{k,n-1}(\nu)+2 P_{k,n}(\nu) - P_{k,n+1}(\nu) = \\
& (P_{k,n}(\nu)-P_{k,n-1}(\nu)) - (P_{k,n+1}(\nu)-P_{k,n}(\nu)) \\
&= m_n(\nu^{k-1})-2 m_n(\nu^k) + m_n(\nu^{k+1}) +
	\sum_a \delta_{\mu_a,n} \delta_{\eta_a,k},
\end{split}
\end{equation*}
valid for $k,n\ge 1$.  In particular the vacancy numbers have the
partial convexity property
\begin{equation} \label{convex}
  P_{k,n}(\nu) \ge 1/2 (P_{k,n-1}(\nu) + P_{k,n+1}(\nu))
  \qquad\text{ if $m_n(\nu^k)>0$.}
\end{equation}
So for 1 it is enough to show that $P_{k,n}(\nu)\ge0$ for sufficiently
large and small $n$.  For small $n$, recall that $P_{k,0}(\nu)=0$.
For large $n$, $Q_n(\nu^k)=|\nu^k|$, so that the vacancy numbers satisfy
\begin{equation} \label{sharp estimate}
\begin{split}
  P_{k,n}(\nu) &= |\nu^{k-1}| - 2 |\nu^k| + |\nu^{k+1}|
  	+ \sum_a \min(\mu_a,n)\delta_{\eta_a,k} \\
  &= (\sum_{j>k-1} \la_j - \sum_{a\ge1} \mu_a (\eta_a-(k-1))_+ \\
  &\,\,\,- 2 (\sum_{j>k} \la_j - \sum_{a\ge1} \mu_a (\eta_a-k)_+ \\
  &\,\,\,(\sum_{j>k+1} \la_j - \sum_{a\ge1} \mu_a (\eta_a-(k+1))_+ 
  	+ \sum_a \min(\mu_a,n)\delta_{\eta_a,k} \\
  &= \la_k - \la_{k+1} + \sum_{a\ge1} \min(0,n-\mu_a) \delta_{\eta_a,k}.
\end{split}
\end{equation}
Thus for large $n$, $P_{k,n}(\nu) \ge \la_k-\la_{k+1}\ge0$, proving 1.

To prove 3, note that
\begin{equation} \label{diff}
  \min(\la_k,n) - \min(\la_{k+1},n) =
  \begin{cases}
    0 & \text{for $n\le \la_{k+1}$} \\
    n - \la_{k+1} & \text{for $\la_{k+1}< n< \la_k$}\\
    \la_k -\la_{k+1} &  \text{for $\la_k \le n$}
  \end{cases}
\end{equation}
For $n\le\la_{k+1}$ there is nothing to show.  Suppose
$n\ge\la_k$.  In light of \eqref{sharp estimate} it suffices to show
that there is no index $a$ such that $\eta_a = k$ and $n<\mu_a$.
Suppose such an $a$ exists.  Then $\mu_a>n\ge\la_k$.  In particular,
the rectangle $R_a$ having $\eta_a=k$ rows and $\mu_a$ columns is not
contained in the shape $\la$.  By Theorem \ref{RC number} it follows that
there are no admissible configurations of type $(\la;R)$, which is a
contradiction.  

Finally suppose $\la_{k+1} < n <\la_k$.  In light of the partial convexity
\eqref{convex} and the boundary conditions $P_{k,\la_{k+1}}(\nu)\ge 0$ and
$P_{k,\la_k}(\nu) \ge \la_k - \la_{k+1}$, it is enough to show that
$m_n(\nu^k)=0$ for $n>\la_{k+1}$.  Since $\nu^k$ is a partition it is
enough to show $\la_{k+1} \ge \nu^k_1$ for all $k$.  Fix $k$.

Observe first that for any $(\la;R)$, if there is an admissible
configuration of type $(\la;R)$ then
$\la_1\ge\mu_a$ for all $a$.  Indeed, if such a configuration exists
then by Theorem \ref{RC number}, $\LRC^R_\la>0$, and by
the LR rule there must be an $R$-LR tableau of shape $\la$.
But for this to happen, $\la$ must contain $R_i$ for all $i$,
hence $\la_1\ge \mu_a$ for all $a$.

Thus it suffices to exhibit a pair $(\lt;\Rt)$ such that
$\lt_1=\la_{k+1}$ and $\Rt$ contains a rectangle having $\nu^k_1$ columns.
Define the partition $\lt$ and sequence of rectangles $\Rt$ by
\begin{equation*}
\begin{split}
  \lt&=(\la_{k+1},\la_{k+2},\dots) \\
\Rt&=\{(\mu_a)^{(\eta_a-k)_+}\} \cup \{(\nu^k_b)\mid 1\le b\le \ell(\nu^k)\}.
\end{split}
\end{equation*}
In other words,
$\lt$ is obtained from $\la$ by removing the first $k$ parts, and
$\Rt$ is obtained from $R$ by removing $k$ rows from each rectangle
and adding a single row for each part of the partition $\nu^k$.
Now it is enough to exhibit an admissible configuration of type
$(\lt;\Rt)$.  Recall that we are assuming that $\nu$ is an admissible
configuration of type $(\la;R)$.  Let $\nt$ be the sequence of partitions
defined by $\nt^p=\nu^{k+p}$ for $p\ge 1$.  Let us check the
condition that $\nt$ is a configuration of type $(\lt;\Rt)$.
For $p\ge1$,
\begin{equation*}
\begin{split}
  |\nt^p| &= |\nu^{k+p}| \\
  &= \sum_{j>k+p} \la_j - \sum_a \mu_a \max(\eta_a-k-p,0) \\
 &= \sum_{j>p} \lt_j - \sum_a \mu_a \max((\eta_a-k)_+-p,0)
-\sum_{1\le b\le \ell(\nu^k)} \nu^k_b \max(1-p,0).
\end{split}
\end{equation*}
To check that $\nt$ is admissible, let $n\ge1$ and $p>1$ first:
\begin{equation*}
\begin{split}
  P_{p,n}(\nt) &= Q_n(\nt^{p-1})-2Q_n(\nt^p)+Q_n(\nt^{p+1})
+\sum_a \min(\mu_a,n) \delta_{(\eta_a-k)_+,p}
+\sum_b \min(\nu^k_b,n) \delta_{1,p} \\
  &= Q_n(\nu^{k+p-1})-2Q_n(\nu^{k+p})+Q_n(\nu^{k+p+1})
+\sum_a \min(\mu_a,n) \delta_{\eta_a,k+p} \\
  &= P_{k+p,n}(\nu) \ge 0
\end{split}
\end{equation*}
by the admissibility of $\nu$.  For $p=1$ and $n\ge1$ we have
\begin{equation*}
\begin{split}
  P_{1,n}(\nt) &= Q_n(\nt^0)-2Q_n(\nt^1)+Q_n(\nt^2)
+\sum_a \min(\mu_a,n) \delta_{(\eta_a-k)_+,1}
+\sum_b \min(\nu^k_b,n) \delta_{1,1} \\
  &= 0-2Q_n(\nu^{k+1})+Q_n(\nu^{k+2})
+\sum_a \min(\mu_a,n) \delta_{\eta_a,k+1} + Q_n(\nu^k) \\
  &= P_{k+1,n}(\nu) \ge 0
\end{split}
\end{equation*}
again by the admissibility of $\nu$.  Thus $\nt$ is
an admissible configuration of type $(\lt;\Rt)$ and the proof is
complete.
\end{proof}

\end{document}